\documentclass[12pt]{amsart}
\usepackage{amssymb, amscd, euscript,hyperref}
\usepackage[matrix,arrow,curve]{xy}
\makeatletter \@addtoreset{equation}{section} \makeatother

\renewcommand\labelenumi{\rm(\roman{enumi})}

\newcommand{\xref}[1]{{\rm~\ref{#1}}}
\newcommand{\stackunder}[2]{\mathrel{\mathop{#2}\limits_{#1}}}

\newcommand{\type}[3]{$\mathbf{(#1{,}#2)^{#3}}$}
\newcommand{\types}[1]{$\mathbf{(2{,}1)^{-}_{#1}}$}
\newcommand{\typee}[1]{$\mathbf{(2{,}0)^{-}_{#1}}$}

\newcommand{\mt}[1]{\operatorname{#1}}
\newcommand{\Sing}{\operatorname{Sing}}
\newcommand{\Supp}{\operatorname{Supp}}

\newcommand{\Bs}{\operatorname{Bs}}
\newcommand{\Pic}{\operatorname{Pic}}
\newcommand{\red}{\operatorname{red}}

\newcommand{\Mori}{\overline{\operatorname{NE}}}
\newcommand{\gen}{\operatorname{gen}}
\newcommand{\Ex}{\operatorname{Ex}}
\newcommand{\mult}{\operatorname{mult}}

\renewcommand{\emptyset}{\varnothing}

\newcommand{\CC}{{\mathbb C}}
\newcommand{\ZZ}{{\mathbb Z}}
\newcommand{\QQ}{{\mathbb Q}}
\newcommand{\PP}{{\mathbb P}}
\newcommand{\FF}{{\mathbb F}}
\newcommand{\NN}{{\mathbb N}}
\newcommand{\RR}{{\mathbb R}}
\newcommand{\OOO}{{\mathcal O}}
\newcommand{\HHH}{{\EuScript{H}}}
\newcommand{\LLL}{{\EuScript{L}}}
\newcommand{\EEE}{{\EuScript{E}}}
\newcommand{\MMM}{{\EuScript{M}}}
\newcommand{\muu}{\mbox{\boldmath $\mu$}}
\newcommand{\ov}[1]{\overline{#1}}
\newcommand{\cone}{\mathbb W}

\newtheorem{theorem}[equation]{Theorem}
\newtheorem{claim}[equation]{Claim}

\newtheorem{proposition}[equation]{Proposition}
\newtheorem{proposition-definition}[equation]{Proposition-definition}
\newtheorem{lemma}[equation]{Lemma}
\newtheorem{sublemma}[equation]{Sublemma}
\newtheorem{corollary}[equation]{Corollary}
\newtheorem{conjecture}[equation]{Conjecture}

\theoremstyle{definition}
\newtheorem{definition}[equation]{Definition}
\newtheorem{example}[equation]{Example}
\newtheorem{pusto}[equation]{}
\newtheorem{setup}[equation]{Assumption}
\newtheorem{remark}[equation]{Remark}

\title{On the degree of Fano threefolds with canonical
Gorenstein singularities}
\author{Yuri G. Prokhorov}
\thanks{The work was
partially supported by grants RFFI 02-01-00441, NS-489.2003.1,
NS-1910.2003.1, and OM RAN}
\address{Department of Algebra, Faculty of Mathematics, Moscow State
Lomonosov University, Moscow 117234, Russia}
\email{prokhoro@mech.math.msu.su}

\dedicatory{To the memory of Andrei Nikolaevich Tyurin}

\begin{document}
\begin{abstract}
We consider Fano threefolds $V$ with canonical Gorenstein
singularities. A sharp bound $-K_V^3\le 72$ of the degree is
proved.
\end{abstract}

\maketitle

\section{Introduction}
In this paper we study Fano threefolds with canonical Gorenstein
singularities. Such varieties appears naturally in the minimal
model theory that explains the following result due to Alexeev:
\begin{theorem}[\cite{A}]
\label{th-alexeev}
Let $Y$ be a $\QQ$-Fano threefold with $\QQ$-factorial terminal
singularities and Picard number $1$. If the anti-canonical model
$\Phi_{|-K_Y|}(Y)$ is three-dimensional, then $Y$ is birationally
equivalent to a Fano threefold $V$ with canonical Gorenstein
singularities and base point free linear system $|-K_V|$.
\end{theorem}
G. Fano \cite{Fano}, \cite{Fano1} studied algebraic threefolds
$V\subset\PP^n$ whose general linear sections $V\cap \PP^{n-2}$
are canonical curves. Under some additional assumptions, these
varieties also are Fano threefolds with canonical Gorenstein
singularities (see \cite{CM}).

By the main result of \cite{Bor} the degree of Fano threefolds
with canonical Gorenstein singularities is bounded. An explicit
bound $-K^3\le 184$ was obtained by Cheltsov \cite{Ch}. However it
is very far from being sharp:

\begin{conjecture}[Fano, Iskovskikh]
Let $V$ be a Fano threefold with canonical Gorenstein
singularities. Then $-K_V^3\le 72$.
\end{conjecture}

A sharp bound of the degree was known under additional
restrictions to singularities:

\begin{theorem}[\cite{Namikawa}]
\label{th-Namik}
Let $V$ be a Fano threefold with terminal Gorenstein
singularities. Then $X$ can be deformed to a nonsingular Fano
threefold. In particular, $-K_{V}^3\le 64$.
\end{theorem}

Moreover, Mukai's vector bundle method can be applied in this case
to obtain a classification of Fano threefolds with terminal
Gorenstein singularities \cite{mukai}. In the canonical case the
situation is more complicated: the assertion of Theorem
\xref{th-Namik} is no longer true that shows the following
example.

\begin{example}[cf. {\cite{Fano}}, {\cite[Ch. 4, Remark 4.2]{Isk0}}]
\label{ex-main}
Weighted projective spaces $\PP(3,1,1,1)$ and $\PP(6,4,1,1)$ are
Fano threefolds with canonical Gorenstein singularities. Here
$-K^3=72$, so they cannot be deformed to smooth ones.
\end{example}

The main result of this paper is the following theorem.

\begin{theorem}
\label{th-main}
Let $V$ be a Fano threefold with canonical Gorenstein
singularities. Then $-K_{V}^3\le 72$. Moreover, if the equality
$-K_{V}^3= 72$ holds, then $V$ is isomorphic to a weighted
projective space in Example~\xref{ex-main}.
\end{theorem}

Thus the present paper completely proves Fano-Iskovskikh
Conjecture. The following is an immediate consequence of our
results.

\begin{corollary}
\label{cor-x}
Let $V\subset\PP^n$ be a normal projective threefold such that the
following condition holds:
\begin{enumerate}
\item[$(\dag)$]
a general hyperplane section $V\cap \PP^{n-1}$ is a K3 surface
with at worst Du Val singularities.
\end{enumerate}
If $\deg V >72$, then $V$ is a cone.
\end{corollary}
Indeed, according to \cite{CM} (see also \cite{Ch1}) the variety
$V$ is Gorenstein and the anti-canonical class $-K_V$ is the class
of hyperplane section. If the singularities of $V$ are rational,
then they are canonical and by Theorem~\xref{th-main} \ $\deg
V=-K_V^3\le 72$. If the locus of non-rational singularities is
non-empty, then it is zero-dimensional and by the main result of
\cite{Is1} we obtain that $V$ is a cone.

Similarly we have
\begin{corollary}
\label{cor-xx}
Let $U\subset\PP^n$ be a normal projective threefold such that the
following condition holds:
\begin{enumerate}
\item[$(\dag\dag)$]
a general hyperplane section $U\cap \PP^{n-1}$ is an Enriques
surface with at worst Du Val singularities.
\end{enumerate}
If $\deg U >36$, then $U$ is a cone.
\end{corollary}
Indeed, assume that $U$ is not a cone. Then according to
\cite{Ch1} (cf. \cite{CM1}), the variety $U$ is $\QQ$-Gorenstein,
has only canonical singularities, and its anti-canonical Weil
divisor $-K_U$ is $\QQ$-linearly equivalent to the class of
hyperplane section $H$. This means that $n(K_U+H)\sim 0$ for some
$n\in \NN$. The divisor $K_U+H$ defines an $n$-sheeted covering
$\pi\colon V\to U$ which is \'etale over the smooth locus of $U$.
Then the variety $V$ satisfies the conditions of Theorem
\xref{th-main} and therefore, $-K_{V'}^3\le 72$. Hence, $\deg
U=-K_U^3=-K_{V}^3/n\le 36$.

Note that in contrast with Corollary~\xref{cor-x}, the bound in
Corollary~\xref{cor-xx} is not sharp. For example it is easy to
show that an involution can not act on varieties from Example
\xref{ex-main} so that the quotient has only canonical
singularities and the action is free in codimension two.
Therefore, $\deg U \neq 36$.

We give one more consequence of our theorem.
\begin{corollary}
Let $X$ be a Fano threefold with canonical \textup(not necessarily
Gorenstein\textup) singularities. Assume that
\begin{enumerate}
\item
the image $\Phi_{|-K_X|}(X)$ of the map given by the
anti-canonical linear system is three-dimensional;
\item
a general element $F\in |-K_X|$ is irreducible and has at worst Du
Val singularities.
\end{enumerate}
Then $\dim |-K_X|\le 38$ and the inequality is sharp.
\end{corollary}
The proof of the corollary is similar to that of Corollary 4.6 in
\cite{A}: consider a $K_X+|-K_X|$-crepant modification $f\colon
(Y,\LLL_Y)\to (X,|-K_X|)$ of the pair $(X,|-K_X|)$. Then the
linear system $\LLL_Y\subset |-K_Y|$ is base point free and
defines a morphism $g\colon Y\to \PP^N$ whose Stein factorization
$Y\to V\to g(Y)$ gives us is a Fano threefold $V$ with canonical
Gorenstein singularities. By our theorem $\dim |-K_X|\le \dim
|\LLL_Y|\le \dim |-K_V|\le 38$.

The last corollary is an argument justifying the conjecture that
the bound in Theorem~\xref{th-main} holds for all Fano threefolds
with canonical singularities (H. Takagi).

In conclusion, we note that in the case, when the variety $V$ has
only cDV singularities, the bound $-K_V^3\le 72$ is not sharp. We
expect in this case the better bound $-K_V^3\le 64$. The proof can
be obtained by a small modification of our method. On the other
hand, the following conjecture looks very realistic:
\begin{conjecture}[cf. Theorem {\xref{th-Namik}}]
Let $V$ Fano threefold with only cDV singularities. Then $V$ has a
smoothing by a small deformations. In particular, $-K_V^3\le 64$.
\end{conjecture}

\begin{remark}
Rationality questions for Fano threefolds satisfying conditions of
Theorem~\xref{th-main} were discussed in \cite{P-F}.
\end{remark}

\begin{pusto}
\label{intro-del-Pezzo}
We illustrate the proof of Theorem~\xref{th-main} on its
two-dimensional analogue. More precisely, we reprove a well-known
fact: the degree of a del Pezzo surface $V$ with Du Val
singularities is bounded by $8$ (if $V$ is singular). Let $V$ be a
del Pezzo surface with Du Val singularities of degree $K_V^2=d$.
Assume that $V$ is singular. If $d\ge 3$, then the anti-canonical
linear system defines an embedding $V \hookrightarrow \PP^d$.
Consider the linear system $\HHH\subset |-K_V|$ of hyperplane
sections passing through a singular point $P\in V$. Let
$\phi\colon W\to V$ be the minimal resolution. Then
$\phi^*(K_V+\HHH)=K_W+\HHH_W+B$, where $\HHH_W$ is the proper
transform of $\HHH$ and $B$ is a nonzero effective divisor. Run
the Minimal Model Program with respect to $K_W+\HHH_W$. Since the
linear system $\HHH_W$ has no fixed components, we do not leave
the category of nonsingular surfaces. At the end we get a pair
$(X,\HHH_X)$ having a nonbirational $K_X+\HHH_X$-negative
contraction. There are two possibilities.

a) $X\simeq \PP^2$. Then $d-1=\dim \HHH= \dim \HHH_X\le 5$.

b) $X\simeq\FF_n$. Then $\HHH_X$ is a linear system of sections of
the projection $\FF_n\to\PP^1$.

Note that the birational transform of the linear system $|-K_V|$
is contained in $|-K_X|$. In particular, $|-K_X|$ has no fixed
components. In case b) this is possible only if $n\le 2$ and then
$d=\dim |-K_V|\le \dim |-K_X|\le 8$.
\end{pusto}

Our proof in the three-dimensional case follows the one outlined
above. However it is much more complicated because of numerous
technical details. Let us explain three-dimensional case in
details.
\begin{definition}
A projective algebraic variety $W$ is called a \emph{weak Fano}
variety, if its anti-canonical divisor is a nef and big
$\QQ$-Cartier divisor.
\end{definition}

Let $V$ be a Fano threefold with canonical Gorenstein
singularities. Consider its terminal $\QQ$-factorial modification
$\phi\colon W\to V$ (see Proposition~\xref{prop-def-term-modif}).
Here $K_W=\phi^*K_V$ and $W$ has only terminal $\QQ$-factorial
Gorenstein singularities. According to \cite[Lemma 5.1]{Kaw-crep}
any Weil divisor on $W$ is Cartier, i.e., $W$ is factorial.
Conversely, for any weak Fano threefold $W$ with terminal
factorial singularities, its multiple anti-canonical image
$V=\Phi_{|-nK_W|}(W)$, for some $n\in \NN$, is a Fano threefold
with canonical Gorenstein singularities. Moreover, $K_W=\Phi^*K_V$
and $-K_V^3=-K_W^3$. Thus the inequality in Theorem~\xref{th-main}
is equivalent to the same inequality for the degree of weak Fano
threefolds $W$ with terminal factorial singularities. On such
varieties $K$-negative extremal rays are completely classified
\cite{Cu}. Simple analysis of the structure of extremal rays shows
that, for $W$, having the maximal degree $\ge 72$, the
corresponding variety $V$ either is singular along a line or
contains a plane. Further, as in our two-dimensional analog the
bound of the degree of $V$ is obtained by detailed analysis of the
linear system of hyperplane sections passing through this line or
plane.

The work was conceived during the Fano Conference (Turin, 2002)
and basically completed at the Kyoto Research Institute for
Mathematical Sciences (RIMS) in 2002-2003. The author is grateful
to organizers of the conference for the invitation and RIMS for the
support and very nice working environment. The author also would
like to thank H. Takagi, D. Orlov, A. Kuznetsov, I. Cheltsov, P. Jahnke, and I.
Radloff for some useful discussions, as well as, 
Haidong Liu and Wenyou Li for pointing out a gap in the journal version 
of this paper.

\section{Notation, conventions and preliminary results}
\label{sect-predv-s2}
\begin{pusto}
In this paper we work over $\CC$, the complex number field. All
Fano varieties are usually supposed to be three-dimensional.
However Lemmas~\xref{lemma-pairs-gen-0} and
\xref{lemma-pairs-gen-1} hold in arbitrary dimension modulo the
Minimal Model Program (MMP). Always when we say that a variety
has, say, canonical singularities, it means that singularities are
not worse than that. Usually we do not distinguish between Cartier
divisors and corresponding invertible sheaves.
\end{pusto}

\subsection*{Notation}
By $\FF_n$ we denote a rational scroll (Hirzebruch surface) and by
$\Sigma=\Sigma_n$ and $l=l_n$ we denote its minimal section and
fiber, respectively. For $n\ge 1$, there exists a contraction
$\FF_n\to \cone_n$ of the minimal section $\Sigma$, where
$\cone_n$ is the cone over a rational normal curve of degree $n$
in $\PP^n$. The vertex of this cone will be denoted by $O=O_n$.

\begin{pusto}
\label{not-Fano}
Everywhere below we assume that $V$ is a Fano threefold with
canonical Gorenstein singularities. Using the Riemann-Roch formula
and Kawamata-Viehweg Vanishing Theorem, it is easy to find the
dimension of the anti-canonical linear system (see \cite{reid0}):
\begin{equation}
\label{lemma-g-dim}
\dim |-K_V|=-\frac12K_V^3+2.
\end{equation}
Denote $g:=-\frac12K_V^3+1$. This number is called the
\emph{genus} of $V$. Thus,
\[
-K_V^3=2g-2\quad \text{and}\quad \dim |-K_V|=g+1.
\]
\end{pusto}

\begin{theorem}[\cite{reid0}]
\label{th-k3}
A general member $L\in |-K_V|$ has only Du Val singularities.
\end{theorem}

The following proposition is an easy consequence of Theorem
\xref{th-k3} and corresponding facts on K3 surfaces \cite{SD}.
\begin{proposition}[{\cite[Ch. 1, Proposition 6.1]{Isk0}}, \cite{Shin}]
\label{2.4.1}
In the above notation, assume that $\Bs |-K_V|\neq \emptyset$.
Then one of the following holds:
\begin{enumerate}
\item
the scheme $\Bs |-K_V|=C$ is a \textup(reduced\textup) nonsingular
rational curve and a general divisor $L\in |-K_V|$ is also
nonsingular along $C$; in particular, $C \cap \Sing V=\emptyset$;
\item
$\Bs |-K_V|=\{P\}$ is a point, a general divisor $L\in |-K_V|$ has
at $P$ an ordinary double singularity and $P\in \Sing V$
\textup(in this case $-K_V^3=2$\textup).
\end{enumerate}
\end{proposition}

\begin{proposition}[\cite{Ch}, cf. {\cite[Ch. 1, \S 6]{Isk0}},
see also \cite{JR}] If $\Bs |-K_V|\neq \emptyset$, then $-K_V^3\le
46$.
\end{proposition}

Similar to the nonsingular case (see \cite{Isk0}) one can easily
prove the following.

\begin{proposition}
\label{prop-hyperell}
In the above notation, let $\Phi=\Phi_{|-K_V|}\colon V
\dashrightarrow \ov V\subset \PP^{g+1}$ be the anti-canonical map,
where $\ov V=\Phi(V)$. Then $\dim \ov V\ge 2$ and $\dim \ov V=2$
if and only if $\Bs |-K_V|\neq \emptyset$. If $\Bs
|-K_V|=\emptyset$, then one of the following holds:
\begin{enumerate}
\item
$\Phi\colon V\to \ov V$ is double cover, in this case $\ov
V\subset \PP^{g+1}$ is a variety of degree $g-1$;
\item
$\Phi\colon V\to \ov V$ is an isomorphism.
\end{enumerate}
\end{proposition}

\begin{proposition}[\cite{Ch}, cf. {\cite[Ch. 2, \S 2]{Isk0}}]
In case \textup{(i)} of Proposition~\xref{prop-hyperell} the
following inequality holds: $-K_V^3\le 40$.
\end{proposition}

\section{Some facts from the minimal model theory}
Basic definitions and facts from the minimal model theory can be
found in numerous surveys and textbooks (see, e.g., \cite{KM1}).
Here we recall only these facts which are not included in the most
of surveys.

\subsection*{Singularities of linear systems (see \cite{A})}
Let $X$ be a normal variety and let $\HHH$ be a linear system (of
Weil divisors) without fixed components on $X$. For any birational
map $\psi \colon X \dashrightarrow Y$, we denote the proper
(birational) transform $\psi_*\HHH$ of $\HHH$ by $\HHH_Y$.
Sometimes, if it does not cause confusion, we will write $\HHH$
instead of $\HHH_Y$. Assume that $K_X+H$ is a $\QQ$-Cartier
divisor for $H\in \HHH$. For any good resolution $f\colon Y\to X$
of singularities of the pair $(X,\HHH)$ we can write
\[
K_Y+\HHH_Y=f^*(K_X+\HHH)+\sum_E a(E,\HHH)E,
\]
where $E$ runs through all exceptional divisors and $a(E,\HHH)\in
\QQ$. Here and below in numerical formulas we will use notation
$\HHH$, $\HHH_Y$ and so on, instead of general members of these
linear systems. We say that a pair $(X,\HHH)$ has \emph{canonical}
singularities (or simply a pair $(X,\HHH)$ is \emph{canonical}),
if $a(E,\HHH)\ge 0$ for all $E$. If a general member of $\HHH$ is
irreducible, then a pair $(X,\HHH)$ is canonical if and only if so
is the pair $(X,H)$, where $H\in \HHH$ is a general divisor
\cite[Th. 4.8]{pairs}. The Log Minimal Model Program works in the
category of three-dimensional canonical $\QQ$-factorial pairs
\cite{A}.

\begin{lemma}
\label{lemma-crep-otobr-1}
Let $(X,\LLL)$ be a canonical pair and let $X\dashrightarrow Y$ be
a birational map such that the inverse map does not contract any
divisors \textup(birational $1$-contraction in Shokurov's
terminology\textup). Assume that $K_X+\LLL\equiv 0$. Then the pair
$(Y,\LLL_Y)$ is also canonical and $K_{Y}+\LLL_Y\equiv 0$.
\end{lemma}

\begin{proof}
Consider a ``Hironaka hut''
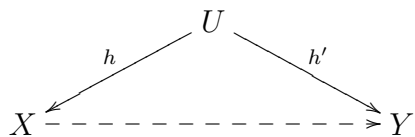
\[
\xymatrix{
&&U\ar[dll]_{h}\ar[drr]^{h'}&&
\\
X\ar@{-->}[rrrr]&&&&Y
}
\]
Let $\{E_i\}$ be the set of all $h'$-exceptional divisors. Write
\begin{eqnarray*}
K_U+\LLL_U&=&h^*(K_X+\LLL_X)+\sum a_iE_i,\quad a_i\ge 0,
\\
K_U+\LLL_U&=&h^{\prime *}(K_{Y}+\LLL_{Y})+\sum b_i E_i.
\end{eqnarray*}
Then $\sum b_i E_i\stackunder{Y}{\equiv} \sum a_iE_i$. This gives
us $a_i=b_i$ (see \cite[1.1]{Sh}).
\end{proof}

\begin{definition}
Let $(X,\LLL)$ and $(X',\LLL')$ be log pairs. A birational
morphism $f\colon (X,\LLL)\to (X',\LLL')$ is said to be
\emph{$K+\LLL$-crepant}, if
\[
K_X+\LLL=f^*(K_{X'}+\LLL'),\qquad \LLL'=f_*\LLL.
\]
A birational map $(X,\LLL)\dashrightarrow (X',\LLL')$ is said to
be \emph{crepant}, if there exists a normal variety $U$ and a
commutative triangle of $K+\LLL$-crepant birational morphisms
\[
\xymatrix{
&&U\ar[dll]_{h}\ar[drr]^{h'}&&
\\
X\ar@{-->}[rrrr]&&&&X'
}
\]
\end{definition}

\begin{lemma}[{\cite[Lemma 3.10]{pairs}}]
\label{lemma-crep-otobr-2}
Let $(X,\LLL)$ be a canonical pair and let $X\dashrightarrow Y$ be
a crepant birational map. Then the pair $(Y,\LLL_Y)$ is also
canonical.
\end{lemma}

\begin{lemma}
\label{lemma-sohranyaetsya}
Let $(W, \HHH)$ be a canonical pair. Assume that
\begin{enumerate}
\item
$X$ has only $\QQ$-factorial terminal singularities;
\item
the linear system $\HHH$ is nef;
\item
$\HHH$ consists of Cartier divisors.
\end{enumerate}
Then the steps of the $K+\HHH$-Minimal Model Program preserve
properties \textup{(i)-(iii)}.
\end{lemma}
\begin{proof}
Let $\varphi\colon W\to W'$ be an extremal $K_W+\HHH$-negative
contraction and let (as usual) $\HHH'=\varphi_*\HHH$. Since the
linear system $\HHH$ is nef, this contraction is also
$K_W$-negative. Therefore the property (i) is kept under
divisorial contractions and flips.

Assume that the contraction $\varphi$ is divisorial. Then the
variety $W'$ is terminal. Let $C$ be a contractible curve. If
$\HHH\cdot C=0$, then by properties of extremal contractions,
$\HHH=\varphi^*\HHH'$ and $\HHH'$ is a linear system of nef
Cartier divisors on $W'$. Assume that $\HHH\cdot C\ge 1$. Then
$-K_W\cdot C>1$. In particular, $\varphi$ contracts a divisor $E$
to a point $P'\in W'$ (see \cite[2.3]{Mori1}). So the linear
system $\HHH'$ is nef. Indeed, if $R'\subset W'$ is an arbitrary
curve and $R\subset W$ is its proper transform, then $\HHH'\cdot
R'\ge \HHH\cdot R\ge 0$. For a general divisor $H\in \HHH$, the
restriction $\varphi|_H\colon H\to H'=\varphi(H)$ is a birational
contraction of normal surfaces with Du Val singularities (because
pairs $(W,\HHH)$ and $(W',\HHH')$ are canonical). By the
Adjunction Formula the divisor $-K_H=-(K_W+H)$ is
$\varphi|_H$-ample. Hence discrepancies of all exceptional over
$P'\in H'$ divisors are strictly positive \cite[Lemma 3.38]{KM1}.
This means that the point $P'\in H'$ is nonsingular. Therefore so
is the point $P'\in W'$ (see \cite[Corollary 3.7]{Sh}) and $H'$ is
a Cartier divisor at $P'$.

Assume that the contraction $\varphi$ is small (flipping).
Consider the flip $W\stackrel{\varphi}{\longrightarrow}
W'\stackrel{\psi}{\longleftarrow} W^{{+}}$. Let $C$ be an
exceptional curve. According to \cite[2.3]{Mori1} we have
$-K_W\cdot C<1$. Hence $\HHH\cdot C=0$. As above, by properties of
extremal contractions, $\HHH=\varphi^*\HHH'$ and $\HHH'$ is a
linear system of nef Cartier divisors on $W'$. Therefore
properties (ii)-(iii) hold for $\HHH^{{+}}=\psi^*\HHH'$.
\end{proof}

\subsection*{Terminal modification}
\begin{proposition-definition}[\cite{reid}, \cite{Mori1},
{\cite[Th. 6.23, 6.25]{KM1}}]
\label{prop-def-term-modif}
Let $X$ be a threefold with only canonical singularities. Then
there exists a threefold $Y$ with only $\QQ$-factorial terminal
singularities and a birational contraction $f\colon Y\to X$ such
that $K_Y=f^*K_X$. Such an $f$ is called a \emph{terminal
$\QQ$-factorial modification} of $X$.
\end{proposition-definition}

\section{Contractions of extremal rays on weak Fano threefolds}
\label{sec-ape}
\label{sec-0pairs}
To work with Fano varieties it is very convenient to use the
following definition which was introduced in the unpublished
preprint of V.V. Shokurov and the author.

\begin{definition}
\label{def:FT}
A \emph{$0$-pair} is a pair $(X,D)$, consisting of a projective
algebraic variety $X$ and a boundary $D$ on $X$ such that
\begin{enumerate}
\item
$(X,D)$ is Kawamata log terminal;
\item
$K_X+D\equiv 0$.
\end{enumerate}
\end{definition}
We say that $0$-pair $(X,D)$ is \emph{generating}, if $X$ is
$\QQ$-factorial and components of $D$ generate the group $N^1(X)$
of $\RR$-divisors modulo numerical equivalence. The following two
lemmas are easy consequences of the log Minimal Model Program.

\begin{lemma}[\cite{PSh}]
\label{lemma-pairs-gen-0}
Let $(X,D)$ be a generating $0$-pair.
\begin{enumerate}
\item
There exists a boundary $\Delta$ such that the pair $(X,\Delta)$
is Kawamata log terminal, the divisor $-(K_X+\Delta)$ ample and
$\Supp(\Delta)=\Supp(D)$.
\item
The Mori cone $\Mori(X)$ is polyhedral and generated by
contractible extremal rays \textup(i.e., rays $R$, for which there
exists a contraction $\varphi_R\colon X\to X'$ in the sense of
Mori \cite{Mori}, a morphism $\varphi_R$ of normal varieties with
connected fibers such that the image of a curve $C$ is a point if
and only if $[C]\in R$\textup).
\item
$\Theta$-Minimal Model Program works with respect to any
\textup(not necessarily effective\textup) divisor $\Theta$.
\end{enumerate}
\end{lemma}

\begin{lemma}[\cite{PSh}]
\label{lemma-pairs-gen-1}
Let $(X,\Delta)$ be a $\QQ$-factorial Kawamata log terminal pair
such that the divisor $-(K_X+\Delta)$ is nef and big. Then there
exists a boundary $D$ such that $(X,D)$ is a generating $0$-pair.
\end{lemma}

\begin{lemma}
\label{lemma-Kaw-Vieh}
Let $(X,D)$ be a generating $0$-pair and let $M$ be an integral
nef Weil divisor on $X$. Then $H^i(X,M)=0$ for all $i>0$.
\end{lemma}

\begin{proof}
By Lemma~\xref{lemma-pairs-gen-0} there is a boundary $\Delta$
such that the pair $(X,\Delta)$ is Kawamata log terminal and the
divisor $-(K_X+\Delta)$ is ample. Now the statement follows by the
Kawamata-Viehweg vanishing theorem.
\end{proof}

Let $R\subset \Mori(X)$ be a (not necessarily $K$-negative)
extremal ray on a normal projective variety $X$. Put
\begin{equation*}
\Ex (R)=\bigcup_{[C]\in R} C.
\end{equation*}
We say that a contractible extremal ray $R$ is of type
\type{n}{m}{} if $\dim \Ex(R)=n$ and $\dim \varphi_R(\Ex(R))=m$.
Further, we distinguish cases \type{n}{m}{-}, \type{n}{m}{+} and
\type{n}{m}{0} according to the sign of $K_X\cdot R$.

Everywhere below we shall assume that $W$ is a weak Fano threefold
with $\rho(W)>1$ having only terminal factorial singularities.
According to Lemma~\xref{lemma-pairs-gen-1} all the extremal rays
on $W$ are contractible. The anti-canonical divisor $-K_W$ defines
a face of the Mori cone $\Mori(W)$ and its contraction $\phi\colon
W\to V$ to a Fano threefold with canonical Gorenstein
singularities.

Consider a birational $K$-negative contraction $f\colon W\to W'$.
According to \cite[(2.3.2)]{Mori1} it is divisorial.

\begin{proposition-definition}
\label{prop-def-types}
Notation as above. Let $S$ be an exceptional divisor. If $f\colon
W\to W'$ is of type \typee{}, then one of the following holds:

{\bf Case \typee{\bullet}.} $W'$ also has only terminal factorial
singularities and is a weak Fano threefold with $-K_{W'}^3\ge
-K_W^3$;

{\bf Case \typee{\circ}.} $S\simeq\PP^2$,
$\OOO_S(S)\simeq\OOO_{\PP^2}(2)$ and $f(S)\in W'$ is a point of
type $\frac12(1,1,1)$.

If $f\colon W\to W'$ is of type \types{}, then $W'$ also has only
terminal factorial singularities and one of the following holds:

{\bf Case \types{\bullet}.} $W'$ is also a weak Fano threefold
with $-K_{W'}^3\ge -K_W^3$;

{\bf Case \types{\circ}.} For $C:=f(S)$ we have $K_{W'}\cdot C>0$
and $C$ is the only curve having negative intersection with
$K_{W'}$. In this case, there exists a small contraction of $C$,
$C\simeq\PP^1$ and $W'$ smooth along $C$. There are two subcases:
\begin{enumerate}
\item[]
\types{\circ0} \qquad $S\simeq \PP^1\times \PP^1$;
\item[]
\types{\circ1} \qquad $S\simeq \FF_1$.
\end{enumerate}
\end{proposition-definition}

\begin{proof}
In case \typee{}, the assertion follows by the classification of
extremal rays on threefolds with terminal factorial singularities
\cite{Cu}.

Consider case \types{}. According to \cite{Cu} the variety $W'$
smooth along $C$ and $f$ is the blow up of $C$. We have
$K_W=f^*K_{W'}+S$ and
\begin{equation*} K_W^3=K_{W'}^3+3f^*K_{W'}\cdot S^2+S^3.
\end{equation*}
If $W'$ is a weak Fano threefold, then
\begin{equation*}
0\le (-K_W)^2\cdot S=2f^*K_{W'}\cdot S^2+S^3,\quad 0\le
(-K_W)\cdot (-f^*K_{W'})\cdot S=f^*K_{W'}\cdot S^2.
\end{equation*}
This gives us $K_W^3\ge K_{W'}^3$

Assume now that $W'$ is not a weak Fano threefold. Then
$K_{W'}\cdot \Gamma>0$ for some irreducible curve $\Gamma$. Let
$\tilde \Gamma\subset W$ be a irreducible curve dominating
$\Gamma$. Since $f^*K_{W'}\cdot \tilde \Gamma>0$ and $K_W\cdot
\tilde \Gamma\le 0$, we have $S\cdot \tilde \Gamma<0$. Therefore,
$\Gamma=C$. By Lemma~\xref{lemma-pairs-gen-1} there exists a
boundary $D$ such that $(W,D)$ is a generating $0$-pair. Let
$D':=f_*D$. Since $K_W+D\sim 0$, we have $K_{W'}+D'\sim 0$ and
$(W',D')$ is also a $0$-pair. Further, $D'\cdot C<0$ and $C$ is
the only irreducible curve with this property. By Lemma
\xref{lemma-pairs-gen-0} the curve $C$ generates an extremal ray
and the corresponding contraction is small. From the
Kawamata-Viehweg Vanishing Theorem one immediately obtains that
$C\simeq \PP^1$ (see \cite[Corollary 1.3]{Mori1}).

Hence, $S\simeq \FF_n$ for some $n\ge 0$. Further,
\begin{equation*}
K_W^2\cdot S=-K_{W'}\cdot C+2-2p_a(C)=2-K_{W'}\cdot C.
\end{equation*}
It is clear that the restriction $-K_W|_S$ is a section of the
fibration $S\to C$. Let $\Sigma$ be a minimal section and let $l$
be a fiber. Then we can write
\begin{equation}
\label{eq-luchi-K-S}
-K_W|_S\sim \Sigma+(n+a)l.
\end{equation}
Since this divisor is nef, $a\ge 0$. Thus,
\begin{equation}
\label{eq-luchi-K-S-1}
0\le -K_{W'}\cdot C=K_W^2\cdot S-2=(-K_W|_S)^2-2=n+2a-2.
\end{equation}
This gives us $n+2a<2$, i.e., $a=0$ and $n\le 1$.
\end{proof}

\begin{corollary}
Notation as above. We have the following table:
\par\medskip\noindent
\begin{center}
\begin{tabular}{|c||c|c|c|c|c|}
\hline & $S$ & $N_{C/W'}$& $K_{W'}\cdot C$ & $K_W\cdot \Sigma$&
$-K_{W'}^3$
\\
 \hline \hline
\types{\circ0}& $\PP^1\times \PP^1$ & $\OOO(-2) \oplus\OOO(-2)$ &
$2$ & $0$& $-K_W^3-2$
\\
\types{\circ1}& $\FF_1$ & $\OOO(-1)\oplus\OOO(-2)$ & $1$ & $0$&
$-K_W^3$
\\
\hline
\end{tabular}
\end{center}
\end{corollary}
\begin{proof}
In our two cases, we have $-K_W|_S\sim \Sigma+nl$, see
\eqref{eq-luchi-K-S}. From~\eqref{eq-luchi-K-S-1} we get
$-K_{W'}\cdot C=n-2$. Hence,
\[
\deg N_{C/W'}=-K_{W'}\cdot C-2=n-4.
\]
Thus,
\begin{multline*}
K_W^3=K_{W'}^3+3f^*K_{W'}\cdot S^2+S^3=K_{W'}^3-3K_{W'}\cdot
C+2+K_{W'}\cdot C=
\\
= K_{W'}^3-2K_{W'}\cdot C+2= K_{W'}^3+2(n-2)+2=K_{W'}^3+2n-2.
\end{multline*}
\end{proof}

\begin{corollary}
\label{cor-line-plane-0}
In case \types{\circ0}, the image $\phi(S)$ of the exceptional
divisor on $V$ is a line \textup(a rational curve $\Gamma$ such
that $-K_V\cdot \Gamma=1$\textup) and the variety $V$ is singular
along $\phi(S)$.

In cases \types{\circ1} and \typee{\circ},\quad $\phi(S)$ is a
plane \textup(a rational surface $\Pi$ such that $K_V^2\cdot
\Pi=1$\textup).
\end{corollary}

Now we consider contractions of type \type{2}{1}{0}.

\begin{lemma}
\label{lemma-2-factorial}
Let $X$ be a threefold with terminal factorial singularities and
let $f \colon X\to X'$ be an extremal contraction of type
\type{2}{1}{0}. Then the variety $X'$ is $2$-factorial \textup(the
last means that for any integral divisor $G$ on $X'$, $2G$ is a
Cartier divisor\textup).
\end{lemma}

\begin{proof}
Consider a general hyperplane section $H'\subset X'$ and let
$H\subset X$ be a proper transform of $H'$. Pick a point $P\in
X'\cap f(E)$. We may assume that $X'$ is a sufficiently small
neighborhood of $P$. Then $P\in H'$ is a Du Val singularity, and
$f|_H\colon H\to H'$ is its minimal resolution. Consider also a
general section $F'\subset X'$, passing through $f(E)$. Write
$f^*F'=F+rE$, where $F$ is the proper transform and $r\in \NN$.
Restricting to $H$ we get $f_H^* F'_{H'}= F_H+rE_H$. Here
$F'_{H'}$ is a general hyperplane section of the singularity $P\in
H'$. Hence $\Gamma:=rE_H$ is the fundamental cycle of a (Du Val)
singularity $P\in H'$. Write $\Gamma=\sum \gamma_i\Gamma_i$. All
the components $\Gamma_j$ are numerically proportional on $X$. In
particular, all the $\Gamma_j\cdot \Gamma$ are non-zero, i.e., the
fundamental cycle is relatively anti-ample. Assume that $P\in H'$
is not a singularity of type $A_1$. Then, by definition of
fundamental cycle, $(\Gamma-\Gamma_j)\cdot \Gamma_j>0$ for all
$j$. So, $0> \Gamma\cdot \Gamma_j>-2$. Hence, $\Gamma\cdot
\Gamma_j=-1$ for all $j$. On the other hand, $\sum
\gamma_j=-\Gamma^2$ is equal to the multiplicity of a rational
singularity $P\in H'$, i.e., $\sum \gamma_j=2$. This is possible
only if $\Gamma=\Gamma_1+\Gamma_2$ and $P\in H'$ is a singularity
of type $A_2$. Thus we obtain two cases:
\begin{enumerate}
\item
$\Gamma=\Gamma_1$ and $P\in H'$ is a singularity of type $A_1$;
\item
$\Gamma=\Gamma_1+\Gamma_2$, $\Gamma_1\equiv \Gamma_2$ and $P\in
H'$ is a singularity of type $A_2$.
\end{enumerate}
Let $F'$ be an integral divisor on $X'$, let $F$ be its proper
transform on $X$ and let $F\cdot \Gamma_i=a$. Consider, for
example, case (ii). Then $(F+arE)\cdot \Gamma_i= a+ a\Gamma\cdot
\Gamma_i =0$. By the Cone Theorem $F+arE=f^*F'$ and $F'$ is a
Cartier divisor. Similarly, in case (i) we get that $2F'$ is a
Cartier divisor.
\end{proof}

\subsection*{Contractions to a curve}
\begin{proposition}
\label{prop-del-Pezzo}
Assume that on $W$ there exists an extremal ray of type
\type{3}{2}{}. Then $-K_W^3\le 54$.
\end{proposition}

\begin{lemma}
\label{lemma-rho=2-div}
Assume that $\rho(W)=2$ and the morphism $\phi\colon W\to V$
contracts a divisor to a curve. Let $S$ be a prime divisor on $W$.
Then
\[
-K_W^3\le \min (54, \ 4K_W^2\cdot S).
\]
\end{lemma}

\begin{proof}
Put $\bar S:=\phi_*S$, \ $d:=-K_W^3$, \ $\delta:=K_W^2\cdot S$ and
assume that $d>4\delta$. Since the variety $V$ is $\QQ$-factorial
and $\rho(V)=1$, we have $d\bar S\equiv -\delta K_V$. On the other
hand, by Lemma~\xref{lemma-2-factorial} the divisor
$G:=2\phi_*\bar S$ is Cartier. Thus, $-K_V\equiv
\frac{d}{2\delta}G$, where $\frac{d}{2\delta}>2$. Let $H$ be an
ample generator of the group $\Pic(V)\simeq\ZZ$. Then $-K_V=rH$,
where $r\ge 3$. A general member $H\in |H|$ is a del Pezzo surface
with Du Val singularities (see \cite{Shin}). By the Adjunction
Formula, $-K_H\equiv (r-1)H|_H$. Hence, $H\simeq \PP^2$ or $H$ is
a quadric in $\PP^3$. It immediately follows that $V$ is
isomorphic to $\PP^3$ or a quadric in $\PP^4$. The first case is
impossible because $V$ is singular. In the second case, we have
$-K_V^3=54$.
\end{proof}
\begin{proof}[Proof of Proposition {\xref{prop-del-Pezzo}}]
A general fiber $W_\eta$ of the corresponding contraction is a
smooth del Pezzo surface. Therefore, $K_W^2\cdot
W_\eta=K_{W_\eta}^2\le 9$. By Lemma~\xref{lemma-rho=2-div} we have
$-K_W^3\le 54$.
\end{proof}

\section{The conic bundle case}
\label{sec-conic}
\begin{pusto}
\label{ass-conic}
In this section, we shall assume that $W$ is a weak Fano threefold
with terminal factorial singularities such that there exists an
extremal contraction $f\colon W\to Z$ from $W$ to a surface (i.e.,
of type \type32{}). The main result of this section is the
following proposition.
\end{pusto}

\begin{proposition}
\label{prop-conic}
\begin{enumerate}
\item \textup(cf. \cite{Mori-Mukai-86}, \cite{P_F0}\textup)
The surface $Z$ is smooth and is a weak del Pezzo surface.
\item
\textup(cf. \cite{Dem}\textup) If $f\colon W\to Z$ is not a
$\PP^1$-bundle, then $-K_W^3\le 54$.
\item
If $f\colon W\to Z$ is a $\PP^1$-bundle, then $-K_W^3\le 64$ with
a unique exception:
\begin{enumerate}
\item[(*)]
$-K_W^3=72$, $Z\simeq \PP^2$ and
$W\simeq\PP(\OOO_{\PP^2}\oplus\OOO_{\PP^2}(3))$ \textup(see
Example~\xref{ex-main}\textup).
\end{enumerate}
\end{enumerate}
\end{proposition}

\begin{proof}
(i) According to \cite{Cu} the surface $Z$ is smooth and $W/Z$ is
a (possibly singular) conic bundle. Let $\EEE:=f_*\OOO_W(-K_W)$.
Then $\EEE$ is a rank $3$ vector bundle. Consider its
projectivization $\PP(\EEE)$ and let $M$ be the tautological
divisor on $\PP(\EEE)$. There exists an embedding
$W\hookrightarrow \PP(\EEE)$ such that $-K_W=M|_W$ and each fiber
$W_z$, $z\in Z$ is a conic in the fiber $\PP(\EEE)_z$ of
$\PP(\EEE)/Z$ (see \cite{Mori}, \cite{Cu}).

Put $L:=f_*K_W^2$. For $m\gg 0$ the divisor $m^2L=f_*(mK_W)^2$ on
$Z$ is moveable and big. We have the standard formula $-4K_Z\equiv
f_*K_W^2+\Delta$ (see, e.g., \cite{Mori-Mukai-86}). Here $\Delta$
is the discriminant of $f$ (a reduced divisor on $Z$). Thus,
\begin{equation}
\label{(*)}
-4K_Z\equiv L+\Delta.
\end{equation}

Assume that there exists an irreducible curve $C\subset Z$ such
that $K_Z\cdot C>0$. It follows from~\eqref{(*)} that
$(4K_Z+\Delta)\cdot C=-L\cdot C\le 0$. Therefore, $\Delta\cdot C<
0$ and $C^2< 0$. Consider two cases:
\subsubsection*{\underline{$C$ is a component of $\Delta$}}
Then
\begin{multline}
\label{eq-Fano-surf-1}
-2\le 2p_a(C)-2=(K_Z+C)\cdot C\le
\\
\le (K_Z+\Delta)\cdot C = -(3K_Z+L)\cdot C\le 0.
\end{multline}
Since $K_Z\cdot C$ is a positive integer, this gives us a
contradiction.

\subsubsection*{\underline{$C$ is not a component of $\Delta$}}
Then
\begin{multline}
\label{eq-Fano-surf-2}
-2\le 2p_a(C)-2=(K_Z+C)\cdot C\le
\\
\le (K_Z+\Delta+C)\cdot C =-(3K_Z+L)\cdot C+C^2< 0.
\end{multline}
As above we get a contradiction. Therefore, the divisor $-K_Z$ is
nef. From~\eqref{(*)} we obtain that $-K_Z$ is big. This proves
(i).

Let us prove (ii). Our proof is completely similar to that of
\cite{Dem}, where the same bound was proved for smooth Fano
threefolds with a conic bundle structure.

The following lemma is well known, however the author could not
find a suitable reference.

\begin{lemma}
\label{lemma-diskr=0}
Let $f\colon X\to Z$ be a \textup(possibly singular\textup) conic
bundle over a smooth surface. Assume that the discriminant curve
$\Delta$ is a tree of rational curves. Then $\Delta=\emptyset$ and
$f$ is a $\PP^1$-bundle \textup(in particular, $X$ is
smooth\textup).
\end{lemma}

\begin{proof}
According to \cite[Theorem 1.13]{Sark1} there exists a standard
model, i.e., the following commutative diagram
\[
\xymatrix{
X\ar[d]_f&X'\ar[d]^{f'}\ar@{-->}[l]
\\
Z& Z'\ar[l]_{\sigma}
}
\]
where $f'$ is a standard conic bundle, $\sigma$ is a composition
of blowups over the finite set $M:=f(\Sing(X))$ and
$X'\dashrightarrow X$ is a birational map which is an isomorphism
over $X\setminus f^{-1}(M)$. The discriminant curve $\Delta'$ of a
standard conic bundle $f'$ is contained in
$\sigma^{-1}(\Delta)\cup \sigma^{-1}(M)$. Therefore it is also a
tree of rational curves. It follows from the Artin-Mumford exact
sequence that $\Delta'=\emptyset$, i.e., $f'$ is a $\PP^1$-bundle.

\begin{lemma}[cf. {\cite{Sark0}}, {\cite[Lemma 4]{Isk1}}]
\label{lemma-surf-(-1)-curve}
Let $f\colon X\to Z$ be a $\PP^1$-bundle over a non-singular
surface $Z$. Assume that there exists $(-1)$-curve $C\subset Z$
and let $\delta\colon Z\to Z'$ be the contraction of $C$.
\begin{enumerate}
\item
The relative Mori cone $\Mori(X/Z')$ is generated by classes of
two curves: the fiber of the projection $f$ and the minimal
section $\Sigma$ of the scroll $f^{-1}(C)$.
\item
Let $f^{-1}(C)\simeq \FF_n$. Then $K_X\cdot \Sigma=n-1$.
\item
If $n=0$, then there exists a commutative diagram
\[
\xymatrix{
X\ar[d]_f\ar[r]^{\sigma}&X'\ar[d]^{f'}
\\
Z\ar[r]^{\delta}& Z'
}
\]
where $\sigma$ is the contraction of $f^{-1}(C)\simeq
\PP^1\times\PP^1$ to the second generator and $f'$ is a
$\PP^1$-bundle.

\item
If $n=1$, then there exists a commutative diagram
\[
\xymatrix{
X\ar[d]_f\ar@{-->}[r]^{\chi}&X^{{+}}\ar[r]^{\sigma}&X'\ar[d]^{f'}
\\
Z\ar[rr]^{\delta}&& Z'
}
\]
where $\chi$ is a flop, $\sigma$ is an extremal divisorial
contraction of a divisor to a nonsingular point and $f'$ is a
$\PP^1$-bundle.

\item
If $n\ge 1$, then there exists a commutative diagram
\[
\xymatrix{
&U\ar[ld]_{\sigma}\ar[rd]^{\sigma'}&
\\
X\ar[rd]^f&&X'\ar[ld]_{f'}
\\
&Z&
}
\]
where $\sigma$ is the blowup of $\Sigma$, $\sigma'$ is the
contraction of the proper transform of $f^{-1}(C)$ and $f'$ is a
$\PP^1$-bundle. In this case $\sigma^{\prime -1}(C)\simeq \FF_{m}$
with $m<n$.
\end{enumerate}
\end{lemma}

\begin{proof}
Put $D:=f^{-1}(C)$. The statement of (i) is obvious because
$\rho(X/Z')=2$. (ii) follows from the equalities
\[
-2=(K_D+\Sigma)\cdot \Sigma= (K_X+D)\cdot \Sigma-n=K_X\cdot
\Sigma+f^*C\cdot \Sigma-n=K_X\cdot \Sigma+C^2-n.
\]
To prove (iii) one has to check the contractibility criterion of
the surface $D\simeq \PP^1\times\PP^1$ to the second system of
generators. Let us prove (iv) (cf. \cite[\S 4.1]{IP}). Since
$K_X\cdot \Sigma=0$, there exists a flop $\chi\colon X
\dashrightarrow X^{{+}}$ with center along $\Sigma$. Here the
anti-canonical divisor $-K_X$ is nef over $Z'$ and the same is
true for $-K_{X^+}$. The Mori cone $\Mori(X^{{+}}/Z')$ (as well as
$\Mori(X/Z')$) is generated by two extremal rays. One of them is
small and generated by the flopped curve $\Sigma^+$. Another one
defines a $K_{X^+}$-negative contraction $\sigma'\colon X^{{+}}\to
X'$ over $Z'$. The variety $X^{+}$ is nonsingular \cite[Th.
6.15]{KM1} and the fiber of the projection $X^{{+}}\to Z'$ is
two-dimensional and contains $D^{{+}}$, the proper transform of
$D$. It is easy to see that $K_{X^{{+}}}^2\cdot D^{{+}}=K_X^2\cdot
D$ (see \cite[\S 4.1]{IP}). Hence,
\[
K_{X^{{+}}}^2\cdot D^{{+}}=(K_X+D)^2\cdot D-2(K_X+D)\cdot
D^2=K_D^2-2K_D\cdot D|_D=4.
\]
According to the classification of extremal rays on nonsingular
threefolds \cite{Mori} the morphism $\sigma'$ can not contract a
divisor to a singular point. Therefore $X'$ is nonsingular. Since
$f'$ is smooth outside of the fiber $f'^{-1}(\delta(C))$, $f'$ is
smooth everywhere \cite[Th. 3.5]{Mori}. Therefore $f'$ is a
$\PP^1$-bundle. If $f'(D^{{+}})$ is a curve, then $\sigma'$ is the
blowup of the fiber of the morphism $f'$. But then
$f'^{-1}(\delta(C))=D^{{+}}\simeq \PP^1\times\PP^1$ and on
$X^{{+}}$ there is no small contractions. A contradiction proves
(iv).

Finally, let us prove (v). Let $E$ be the exceptional divisor of
the blowup $\sigma$, let $D_U$ be the proper transform of $D$ on
$U$, and let $D':=\sigma'(E)$. It is clear that the divisor $D_U$
satisfies the contractibility criterion onto a curve. Thus there
exists a contraction $\sigma'\colon U\to X'$ over $Z$, where the
variety $X'$ is nonsingular. Since the restriction
$\sigma'|_{E}\colon E\to D'$ is an isomorphism, it is sufficient
to show that $E\simeq \FF_{m}$ with $m< n$. For the normal sheaf
we have a decomposition $N_{\Sigma/X}=\OOO_{\PP^1}(a)\oplus
\OOO_{\PP^1}(b)$, $a\ge b$. Then $m=a-b$. On the other hand, we
have an exact sequence
\begin{equation}
\label{eq-surf-norm-sheaf}
\begin{array}{lllllllll}
0& \longrightarrow& N_{\Sigma/D}& \longrightarrow& N_{\Sigma/X}&
\longrightarrow& N_{D/X}|_\Sigma& \longrightarrow& 0.
\\
&&\|&&&&\|
\\
&&\OOO_{\PP^1}(-n)&&&&\OOO_{\PP^1}(-1)
\end{array}
\end{equation}
This gives us
\[
\begin{array}{l}
a+b= \deg N_{\Sigma/D}+ \deg N_{D/X}|_\Sigma=-n-1.
\end{array}
\]
Hence, $n=-a-b-1$. If $-n>b$, then in~\eqref{eq-surf-norm-sheaf}
the projection $\OOO_{\PP^1}(-n)$ to the second summand of
$N_{\Sigma/X}$ is trivial and $-n=a$. Therefore, $b=-1$ and
$m=a-b=1-n<n$. If $-n\le b$, then $a+1\le 0$ and \[ m=a-b\le
-1-b\le -a-b-2=n-1.
\]
\end{proof}

Now we finish the proof of Lemma~\xref{lemma-diskr=0}. According
to Lemma~\xref{lemma-surf-(-1)-curve} we can contract all the
$(-1)$ curves in $\sigma^{-1}(M)$. More precisely, there exists
the following commutative diagram
\[
\xymatrix{
X^{0}\ar[d]^{f^{0}}&X'\ar[d]^{f'}\ar@{-->}[l]
\\
Z& Z'\ar[l]^{\sigma}
}
\]
where $f^{0}$ is a $\PP^1$-bundle and birational map $X'
\dashrightarrow X^{0}$ is an isomorphism over $X\setminus
{f^{0}}^{-1}(M)$. Therefore, conic bundles $X/Z$ and $X^{0}/Z$ are
isomorphic in codimension one. But then they are isomorphic
everywhere, $X$ is smooth and $f$ is a $\PP^1$-bundle.
\end{proof}

Thus we may assume that in case (ii) the discriminant curve is
non-empty.

\begin{lemma}
\begin{enumerate}
\item
$\Delta\sim -3K_Z-c_1(\EEE)$;
\item
$W\sim 2 M+f^*(-K_Z-c_1(\EEE))$;
\item
$-K_W^3= c_1(\EEE)\cdot (-K_Z+c_1(\EEE))-2c_2(\EEE)$.
\end{enumerate}
\end{lemma}
\begin{proof}
The assertion (i) can be found in \cite[\S 1]{Sark1}. For the
proof of (ii) we note that $W\sim 2 M+f^*G$ for some divisor $G$
on $Z$. Since $K_{\PP(\EEE)}= -3M+f^*(c_1(\EEE)+K_Z)$, by the
Adjunction Formula $G\sim -c_1(\EEE)-K_Z$. Finally, (iii) follows
from (ii) and the fact that $-K_W=M|_W$.
\end{proof}

We may assume that $-K_V$ is very ample (see \cite{Ch}). Then the
linear system $|-K_W|=|M||_W$ defines a morphism $W\to \PP^{g+1}$,
which maps fibers $f$ to conics in $\PP^{g+1}$. Let $\PP^2_z$ be a
fiber of $\PP(\EEE)\to Z$ and let $W_z:=\PP^2_z\cap W$. Since the
linear system $|-K_W||_{W_z}$ defines an embedding, so do
$|M||_{\PP^2_z}$. Therefore, the linear system $|M|$ has no base
points. This means that the vector bundle $\EEE$ is generated by
global sections. In particular, the class $c_1(\EEE)=-3K_Z-\Delta$
is nef and $c_2(\EEE)\ge 0$. Thus,
\begin{equation*}
-K_W^3\le (-3K_Z-\Delta)\cdot (-4K_Z-\Delta)=12 K_Z^2+7K_Z\cdot
\Delta+\Delta^2.
\end{equation*}
By Lemma~\xref{lemma-diskr=0} we have $p_a(\Delta)\ge 1$. Now our
assertion is a consequence of the following.

\begin{sublemma}
Let $Z$ be a smooth weak del Pezzo surface and let $\Delta$ be a
reduced curve on $Z$ such that $p_a(\Delta)\ge 1$ and the divisor
$-3K_Z-\Delta$ is nef. Then $12 K_Z^2+7K_Z\cdot \Delta+\Delta^2\le
54$.
\end{sublemma}

\begin{proof}
Let $\sigma \colon Z\to Z'$ be a contraction of a $(-1)$-curve
$L$. Put $\Delta':=\sigma(\Delta)$. Then $K_Z=\sigma^*K_{Z'}+L$
and $\Delta=\sigma^*\Delta'- aL$, where $a\ge 0$. Therefore,
\begin{equation*}
2p_a(\Delta)-2=(K_Z+\Delta)\cdot \Delta=(K_{Z'}+\Delta')\cdot
\Delta' +(a-1)aL^2\le 2p_a(\Delta')-2,
\end{equation*}
\begin{multline*}
12 K_Z^2+7K_Z\cdot
\Delta+\Delta^2=12(\sigma^*K_{Z'}+L)^2+7(\sigma^*K_{Z'}+L) \cdot
(\sigma^*\Delta'- aL)
\\
+(\sigma^*\Delta'- aL)^2=12 K_{Z'}^2-12+7K_{Z'}\cdot\Delta'
+7a+(\Delta')^2-a^2
\\
\le 12 K_{Z'}^2+7K_{Z'}\cdot\Delta' +(\Delta')^2.
\end{multline*}
Thus we may assume that the surface $Z$ contains no $(-1)$-curves,
i.e., $Z$ is isomorphic to $\PP^2$ or $\FF_n$ with $n=0,\,2$.

Consider the case $Z\simeq\PP^2$. Let $d:=\deg \Delta$. Since the
divisor $-3K_Z-\Delta$ is nef, we have $3\le d\le 9$. Hence,
\begin{equation*}
12 K_Z^2+7K_Z\cdot \Delta+\Delta^2= 108-21d+d^2\le 54.
\end{equation*}

Now consider the case $Z\simeq\FF_n$,\ $n=0,\, 2$. Let $\Sigma$ be
the minimal section and let $l$ be a fiber. We can write
$\Delta\sim \alpha \Sigma +\beta l$ for some $\alpha,\beta\in
\ZZ$. Since $p_a(\Delta)\ge 1$,
\begin{multline*}
0\le (K_{\FF_n}+\Delta)\cdot
\Delta=-n\alpha(\alpha-2)+(\alpha-2)\beta+\alpha(\beta-2-n)=
\\
(\alpha-1)(-n\alpha+2\beta-2)-2.
\end{multline*}
In particular, $\alpha\ge 2$. This gives us $-n\alpha+2\beta-2\ge
1$, \ $2\beta \ge 3+n\alpha$. Further,
\begin{multline*}
12 K_Z^2+7K_Z\cdot \Delta+\Delta^2=
96+14n\alpha-14\beta-7n\alpha-14\alpha-n\alpha^2+2\alpha\beta
\\
=(7-\alpha)(n\alpha-2\beta+14)-2.
\end{multline*}
Since the divisor
\[
-3K_Z-\Delta\sim (3-\alpha)\Sigma+(6+3n-\beta)l
\]
is nef, $\alpha\le 6$. Hence,
\[
(7-\alpha)(n\alpha-2\beta+14)-2\le 11(7-\alpha)-2\le 53.
\]
\end{proof}

\begin{pusto}
\label{ass-conic-iii}
Let us prove (iii). Additionally to~\xref{ass-conic} we now assume
that $f\colon W\to Z$ is a $\PP^1$-bundle (in particular, $W$ is
smooth). Then $W=\PP(\EEE)$ for some rank $2$ vector bundle $\EEE$
on $Z$. Let $L$ be the tautological divisor.
\end{pusto}

By the relative Euler exact sequence
\begin{equation}
\label{eq-conic-p1--k}
-K_W = 2L + f^*(-K_Z-c_1).
\end{equation}
The Hirsch formula gives us
\begin{equation}
\label{eq-p1-l2}
L^2\equiv L\cdot f ^*c_1-f^*c_2.
\end{equation}
Combining we obtain $L^3=c_1^2-c_2$ and
\begin{equation}
\label{eq-p1-k3}
-K_W^3= 6K_Z^2+2c_1^2-8c_2.
\end{equation}
By the Riemann-Roch Theorem and the Serre duality
\begin{equation}
\label{eq-p1-RR}
h^0(\EEE)+h^0(\EEE\otimes \det \EEE^*\otimes \omega_Z)\ge
\frac12\bigl(c_1^2-2c_2-K_Z\cdot c_1\bigr)+2.
\end{equation}

\subsection*{Case $Z=\PP^2$}
Denote $G:=f^*h$, where $h$ is a line on $Z=\PP^2$. If $c_1$ even,
we can put $c_1=-2$. Then $L^2\cdot G=-2$ and $-K_W^3=62-8c_2$.
Assuming $-K_W^3\ge 64$, we obtain $c_2\le -1$. From
\eqref{eq-p1-RR} we have $h^0(\EEE)\ge 1-c_2\ge 2$, i.e., $\dim
|L|\ge 1$. Since $L^2\cdot G<0$, the linear system $|L|$ has a
fixed component. On the other hand, $\Pic(W)=\ZZ\cdot
M\oplus\ZZ\cdot G$. Therefore, $|L-G|\neq \emptyset$. But then
\[
0\le (-K_W)\cdot (L-G)\cdot G=(2L+5G)\cdot (L-G)\cdot G=2L^2\cdot
G+3L\cdot G^2=-1,
\]
a contradiction.

Assume that $c_1(\EEE)$ is odd. Then the anti-canonical divisor
$-K_W$ is divisible by $2$:\quad $-K_W=2D$, where $D=L+ \frac12
f^* (-K_Z-c_1(\EEE))$. In this case, a general member $D\in |D|$
is a weak del Pezzo surface with Du Val singularities (see
\cite{Shin}). Then $-K_W^3=8K_D^2$. Assuming $-K_W^3> 64$, we
obtain $K_D^2=9$, $D\simeq\PP^2$ and $D_D=-K_D$. Note that the
restriction $H^0(W,\OOO_W(D))\to H^0(D,\OOO_D(D))$ is surjective.
Therefore, the divisor $D$ is very ample and defines an embedding
$V=V_9\hookrightarrow \PP^{10}$ so that fibers $f$ are mapped to
lines in $\PP^{10}$. Note that $\phi\colon W\to V$ is an extremal
contraction. If $\phi$ is of type \type{1}{0}{0}, then $V$ has
only terminal Gorenstein singularities. In this case, by Theorem
\xref{th-Namik} we have $-K_W^3\le 64$. Assume that $\phi \colon
W\to V$ is divisorial and let $E$ be an exceptional divisor. Then
$V$ is $\QQ$-factorial and $\rho(V)=1$. Since $D\simeq\PP^2$, $D$
does not meet fibers of $\phi$, i.e., $\phi(D)\cap
\phi(E)=\emptyset$. This is possible only if $\phi(E)$ is a point.
Therefore, $V$ is a cone over $\phi(D)$. In other words, $V\simeq
\PP(3,1,1,1)$. We get case (*).

Further, we consider the case when $Z\simeq \FF_n$ with $n=0$ or
$2$. Let $\Sigma$ and $l$ be the minimal section and a fiber of
$\FF_n$ respectively. We may assume that $c_1=c_1(\EEE)=a\Sigma+b
l$ and $c_2=c_2(\EEE)=c$. If both numbers $a$ and $b$ are even, we
can put $c_1=0$. From~\eqref{eq-conic-p1--k} we obtain
\begin{equation*}
-K_W=2L-f^*K_Z=2D.
\end{equation*}
where $D\sim L+f^*(\Sigma+(1+n/2)l)$. As above a general member
$D\in |D|$ is isomorphic to $\PP^2$. On the other hand, there is
no dominant morphisms from $\PP^2$ to $\FF_n$, a contradiction.

\subsection*{Case $Z\simeq \PP^1\times \PP^1$}
In this case we can permute $a$ and $b$.

\subsubsection*{Case: $a$ is odd, $b$ is even}
We may assume that $c_1=-3\Sigma$. From~\eqref{eq-conic-p1--k},
\eqref{eq-p1-l2} and~\eqref{eq-p1-k3} we get
\begin{multline*}
L^2\cdot f^*\Sigma=0, \quad L^2\cdot f^*l=-3,
\\
-K_W=2L+f^*(5\Sigma+2l), \quad -K_W^3=48-8c.
\end{multline*}
Hence, $c\le -2$. By~\eqref{eq-p1-RR} we have
\begin{equation*}
h^0(\EEE)+h^0(\EEE\otimes \OOO(\Sigma-2l))\ge \frac
12(-2c+3K_Z\cdot \Sigma)+2= -c-1>0.
\end{equation*}
Therefore one of the following holds: $|L|\neq \emptyset$ or
$|L+f^*\Sigma-2f^*l|\neq \emptyset$. On the other hand,
\begin{gather*}
0\le -K_W\cdot L\cdot f^*l= (2L+f^*(5\Sigma+2l))\cdot L\cdot
f^*l=-1,
\\
0\le -K_W\cdot (L+f^*\Sigma-2f^*l)\cdot f^*\Sigma <0,
\end{gather*}
a contradiction.

\subsubsection*{Case $c_1=-\Sigma-l$}
As above,
\begin{multline*}
L^2\cdot f^*\Sigma=L^2\cdot f^*l=-1, \quad -K_W=2L+3f^*(\Sigma+l),
\\
-K_W^3=(2L+f^*(3\Sigma+3l))^3= 52-8c.
\end{multline*}
This gives us $c\le -2$. By~\eqref{eq-p1-RR} we have
\begin{equation*}
h^0(\EEE)+h^0(\EEE\otimes \OOO(-\Sigma-l))\ge 1-c>0.
\end{equation*}
Therefore, $|L|\neq \emptyset$. Since $L^2\cdot f^*\Sigma=L^2\cdot
f^*l=-1$, $|L|$ has a fixed component. Write $L=F+|M|$, where $F$
is the fixed part. Then $M=f^*Q$, where $Q\sim \alpha \Sigma+\beta
l$, $\alpha,\, \beta\ge 0$. This gives us
\begin{equation*}
0\le F\cdot (-K_W)\cdot f^* l= (L-f^*(\alpha \Sigma+\beta l))\cdot
(2L+3f^*(\Sigma+l))\cdot f^* l=1-2\alpha
\end{equation*}
Hence, $\alpha=0$. Similarly we get $\beta=0$ and $F=0$, a
contradiction.

\subsection*{Case $Z\simeq \FF_2$}

\subsubsection*{Case $c_1=-\Sigma$}
Then
\begin{multline*}
L^2\cdot f^*\Sigma=2, \quad L^2\cdot f^*l=-1, \quad
-K_W=2L+f^*(3\Sigma+4l),
\\
-K_W^3=(2L+f^*(3\Sigma+4l))^3=44-8c_2.
\end{multline*}
Hence, $c_2\le -3$. By~\eqref{eq-p1-RR} we have
\begin{equation*}
h^0(\EEE)+h^0(\EEE\otimes \OOO(-\Sigma-4l) )\ge 1-c_2\ge 4.
\end{equation*}
Hence, $\dim |L|\ge 2$ (because $\dim |L|>\dim
|L-f^*(\Sigma+4l)|$). Since $L^2\cdot f^*l<0$, the linear system
$|L|$ has a fixed component. Write $|L|=F+|M|$, where $F$ is the
fixed and $|M|$ is the moveable part. Then
\[
-K_W\cdot (F+M)\cdot f^*l=-K_W\cdot L\cdot f^*l=1,
\]
Hence, $-K_W\cdot M\cdot f^*l\le 1$ and $-K_W\cdot F\cdot f^*l\le
1$. Since $-K_W\cdot f^*\Sigma\cdot f^*l=2$, the divisor
$f^*\Sigma$ can not be a fixed component of $|L|$. Therefore, $F$
has exactly one component which must be horizontal and $M\sim
f^*(\alpha\Sigma+\beta l)$ for some $\alpha,\, \beta\ge0$. Then
\[
-K_W\cdot M\cdot f^*l=-K_W\cdot f^*(\alpha\Sigma+\beta l)\cdot
f^*l=2\alpha.
\]
Thus, $\alpha=0$. But then
\[
-K_W\cdot F\cdot f^*\Sigma=-K_W\cdot (L-\beta f^*l)\cdot
f^*\Sigma=2-2\beta.
\]
The only possibility is $\beta=1$. In this case, $\dim |L|=\dim
|M|=\dim |l|=1$, a contradiction.

\subsubsection*{Case $c_1=-l$}
Then
\begin{multline*}
L^2\cdot f^*\Sigma=-1, \quad L^2\cdot f^*l=0,
\\
-K_W=2L+f^*(2\Sigma+5l), \quad
-K_W^3=(2L+f^*(3\Sigma+5l))^3=48-8c_2.
\end{multline*}
Thus, $c_2\le -3$. By~\eqref{eq-p1-RR} we have
\begin{equation*}
h^0(\EEE)+h^0(\EEE\otimes \OOO(-2\Sigma-3l))\ge -c_2+1\ge 4.
\end{equation*}
Hence, $\dim |L|\ge 1$. Since
\begin{equation*}
-K_W\cdot L\cdot f^*\Sigma=(2L+f^*(2\Sigma+5l))\cdot L\cdot
f^*\Sigma=-1<0,
\end{equation*}
$f^*\Sigma$ is a fixed component of $|L|$. Further,
\begin{equation*}
-K_W\cdot (L-f^*\Sigma)\cdot f^*l=(2L+f^*(2\Sigma+5l))\cdot
(L-f^*\Sigma)\cdot f^*l=0.
\end{equation*}
This contradicts the fact that $-K_W$ is nef and big.

\subsubsection*{Case $c_1=-\Sigma-l$}
Then
\begin{multline*}
L^2\cdot f^*\Sigma=1, \quad L^2\cdot f^*l=-1,
\\
-K_W=2L+f^*(3\Sigma+5l), \quad -K_W^3=(2L+f^*(3\Sigma+5l))^3=
48-8c_2.
\end{multline*}
Hence, $c_2\le -3$. By~\eqref{eq-p1-RR} we have
\begin{equation*}
h^0(\EEE)+h^0(\EEE\otimes \OOO(-\Sigma-3l))\ge 1-c_2\ge 4.
\end{equation*}
Thus, $\dim |L|\ge 1$. Since $L^2\cdot f^*l=-1$, the linear system
$|L|$ has a fixed component. Write $|L|=F+|M|$, where $F$ is the
fixed and $|M|$ is the moveable part. Then
\[
-K_W\cdot (F+M)\cdot f^*l=-K_W\cdot L\cdot f^*l=1,
\]
Hence, $-K_W\cdot M\cdot f^*l\le 1$ and $-K_W\cdot F\cdot f^*l\le
1$. Since $-K_W\cdot f^*\Sigma\cdot f^*l=2$, the divisor
$f^*\Sigma$ can not be a fixed component of $|L|$. Therefore, $F$
has exactly one component which must be horizontal and $M\sim
f^*(\alpha\Sigma+\beta l)$ for some $\alpha,\, \beta\ge0$. Then
\[
-K_W\cdot M\cdot f^*l=-K_W\cdot f^*(\alpha\Sigma+\beta l)\cdot
f^*l=2\alpha.
\]
Thus, $\alpha=0$. It follows that
\[
-K_W\cdot F\cdot f^*\Sigma=-K_W\cdot (L-\beta f^*l)\cdot
f^*\Sigma=1-2\beta<0,
\]
a contradiction.

Now we consider cases when $Z$ contains a $(-1)$-curve.

\begin{lemma}
\label{lemma-conic-last}
If $Z$ contains a $(-1)$-curve, then there exists a weak Fano
threefold $W'$ satisfying assumptions of~\xref{ass-conic},
\xref{ass-conic-iii} and such that $-K_{W'}^3\ge -K_W^3+6$.
Moreover, if $-K_{W'}^3= -K_W^3+6$, then the linear system
$|-nK_{W'}|$ on $W'$ can contract only $f'$-vertical divisors.
\end{lemma}

\begin{proof}
Let $C\subset Z$ be a $(-1)$-curve, let $\delta\colon Z\to Z'$ be
its contraction, and let $D:=f^{-1}(C)$. Then $D\simeq \FF_n$ for
some $n$. According to Lemma~\xref{lemma-surf-(-1)-curve}, for the
minimal section $\Sigma$ of the surface $D\simeq \FF_n$ we have
$0\ge K_W\cdot \Sigma=n-1$. Hence, $n\le 1$. If $n=1$, then again
by Lemma~\xref{lemma-surf-(-1)-curve} there exists the following
commutative diagram
\begin{equation*}
\xymatrix{
W\ar[d]^f\ar@{-->}[r]^{\chi}&W^{{+}}\ar[r]^{\sigma}&W'\ar[d]^{f'}
\\
Z\ar[rr]^{\delta}&& Z'
}
\end{equation*}
where $\chi$ is a flop, $\sigma$ is an extremal divisorial
contraction which contracts a divisor to a smooth point and $f'$
is a $\PP^1$-bundle. Then $W^+$ is again a weak Fano threefold,
satisfying assumptions of~\xref{ass-conic} and
\xref{ass-conic-iii}. Since $\sigma$ contacts a divisor to a
point, the same holds for $W'$. It is clear that
$-K_W^3=-K_{W^+}^3=-K_{W'}^3-8$.

The case $n=0$ can be treated by a similar way using (iii) from
Lemma~\xref{lemma-surf-(-1)-curve}. Here $\chi=\mt{id}$ and
$\sigma$ contracts a divisor $S$ to a fiber of the morphism $f'$.
Therefore, $-K_W^3=-K_{W'}^3-6$. Moreover, if $\sigma(S)$ meets a
$K_{W'}$-trivial curve $C'$, then for its proper transform
$C\subset W$ we have
\begin{equation*}
K_{W}\cdot C=\sigma^*K_{W'}\cdot C+S\cdot C>K_{W'}\cdot C'=0,
\end{equation*}
a contradiction.
\end{proof}
Now we finish the proof of Proposition~\xref{prop-conic} in case,
when $Z$ contains a $(-1)$-curve. By inductive hypothesis we may
assume that $-K_{W'}^3\le 72$. If $-K_W^3\ge 66$, then $-K_{W'}^3=
72$ and $W'$ is such as in (*). This contradicts to the second
statement from Lemma~\xref{lemma-conic-last}. This completes the
proof of Proposition~\xref{prop-conic}.
\end{proof}

\section{Construction}
\begin{pusto}
\label{notation_1}
As above, let $W$ be a weak Fano threefold having only terminal
factorial singularities. Let $\phi\colon W\to V$ be a morphism
defined by the linear system $|-nK_W|$ for $n\gg 0$. Assume that
$-K_W^3\ge 72$. We may assume that

\begin{enumerate}
\renewcommand\labelenumi{(\arabic{enumi})}
\item
$W$ has no extremal rays of type \type{3}{1}{} and \type{3}{2}{}
(see Propositions~\xref{prop-del-Pezzo} and~\xref{prop-conic});
\item
$W$ has no extremal rays of type \typee{\bullet} and
\types{\bullet} (see Proposition~\xref{prop-def-types}).
\end{enumerate}
Since $W$ is a weak Fano, there exists at least one $K$-negative
extremal ray $R$ on $W$. By the above, $R$ is of types
\types{\circ1}, \types{\circ0} or \typee{\circ}.
\end{pusto}

\begin{pusto}
\label{cor-line-plane}
Taking into account Corollary~\xref{cor-line-plane-0} we get that
for $V$ at least one of the following holds:
\begin{enumerate}
\item[(A)]
$V$ is singular along a line $\Gamma$, or
\item[(B)]
$V$ contains a plane $\Pi$.
\end{enumerate}
For further inquiry it is very convenient also distinguish the
following case:
\begin{enumerate}
\item[(0)]
the variety $V$ has at least one singular non-cDV point $P$.
\end{enumerate}
\end{pusto}

Recall that a three-dimensional singularity is said to be
\emph{cDV}, if it is hypersurface and locally up to analytic
coordinate change is given by an equation
$\phi(x,y,z)+t\psi(x,y,z,t)=0$, where $\phi(x,y,z)=0$ is an
equation of a (two-dimensional) Du Val singularity \cite{reid}.

Now our theorem is a consequence of the following.
\begin{proposition}
\label{propo-last-end}
Let $V$ be a Fano threefold having only canonical Gorenstein
singularities. Assume that $V$ satisfies conditions \textup{(0)},
\textup{(A)}, or \textup{(B)} of~\xref{cor-line-plane}. Then
$-K_V^3\le 72$. Moreover, if the equality $-K_{V}^3= 72$ holds,
then $V$ is isomorphic to $\PP(3,1,1,1)$ or $\PP(6,4,1,1)$.
\end{proposition}

If $-K_V^3>54$, then according to \cite{Ch} (see \S
\xref{sect-predv-s2}) the anti-canonical linear system $|-K_V|$
defines an embedding $V\hookrightarrow \PP^{g+1}$. Moreover, its
image $V_{2g-2}\subset \PP^{g+1}$ is an intersection of quadrics
(\cite[Lemma 3]{Ch}). Put $\LLL:=|-K_V|$. Then $\dim \LLL=g+1>38$.
Consider the following linear subsystem $\HHH\subset \LLL$:
\begin{enumerate}
\item[]
in case (0):\qquad $\HHH:=\{ H\in \LLL \mid H\ni P\}$;

\item[]
in case (A):\qquad $\HHH:=\{ H\in \LLL \mid H\supset \Gamma\}$;

\item[]
in case (B):\qquad $\HHH:=\{ H \mid H+\Pi\in \LLL\}$.
\end{enumerate}
It is clear that
\[
\dim\HHH=
\begin{cases}
\dim \LLL-1&\text{in case (0),}
\\
\dim \LLL-2&\text{in case (A),}
\\
\dim \LLL-3&\text{in case (B).}
\end{cases}
\]

Now let $\LLL_W$ and $\HHH_W$ be proper transforms of $\LLL$ and
$\HHH$ respectively.

In case (0), according to \cite{reid} there exists at least one
exceptional divisor $B_i$ with center at $P$ and discrepancy
$a(B_i)=0$. Write $\phi^*\HHH_W=\HHH+B$, where $B=\sum b_iB_i$ is
an (integral) nonzero effective exceptional divisor over $P$.
Thus,
\begin{equation*}
K_W+\HHH_W+B\sim 0.
\end{equation*}

Similarly, in all cases we have
\begin{equation}
\label{eq-HB}
\begin{array}{lll}
K_W+\HHH_W+B&=&\phi^*(K_V+\HHH)\sim 0,\quad \text{(cases (0) and
(A))},
\\
K_W+\HHH_W+B&=&\phi^*(K_V+\HHH+\Pi)\sim 0,\quad \text{(case (B))},
\end{array}
\end{equation}
where $B$ is an integral effective non-zero divisor.

\begin{lemma}
\label{lemma-predv-image}
The image of the variety $V$ under the map $\Phi_{\HHH}$ given by
the linear system $\HHH$ is of dimension $3$.
\end{lemma}
\begin{proof}
Note that in case (0), the map $\Phi_{\HHH}$ is nothing but the
projection from $P$. Assume that $\dim \Phi_{\HHH}(V)\le 2$. Then
$V$ is a cone over $\Phi_{\HHH}(V)$ with the vertex at $P$. But in
this case the singularity $P\in X$ can not be canonical (see,
e.g., \cite[2.14]{reid}).

Consider, for example, case (A) (case (B) is considered in a
similar way). As above we note that $\Phi_{\HHH}$ is the
projection from the line $\Gamma$. If $\Phi_{\HHH}(V)$ is a curve
$C$, then the variety $V$ is a cone over $C$ with the vertex at
$\Gamma$. But then a general member $L\in \LLL$ is also a cone
with the vertex at $L\cap \Gamma$. This contradicts the fact that
$L$ is a K3 surface with Du Val singularities. Let now
$\Phi_{\HHH}(V)$ be a surface. Then a general fiber of the map
$\Phi_{\HHH}\colon V \dashrightarrow \Phi_{\HHH}(V)$ is
one-dimensional. On the other hand, fibers of this map are cut out
on $V$ by planes $\Lambda$ passing through $\Gamma$. Every such a
plane cut out on $V$ a scheme which is an intersection of
quadrics. Therefore, $\Lambda\cap V=\Gamma+\Gamma_{\Lambda}$,
where $\Gamma_{\Lambda}$ is a line and fibers of $\Phi_{\HHH}$ are
such lines $\Gamma_{\Lambda}$. Thus there is a line passing
through a general point of $V$. In this case, it is easy to get
the bound $-K_V^3\le 46$ (see \cite[Lemma 5]{Ch}).
\end{proof}

From Theorem~\xref{th-k3} and Inversion of Adjunction \cite[3.3,
9.5]{Sh} we immediately obtain the following.
\begin{corollary}
\label{cor-k3}
The pair $(V,|-K_V|)$ has only canonical singularities.
\end{corollary}

We can write
\[
K_W+\LLL_W=\phi^*(K_V+\LLL)\sim 0.
\]
Hence the pair $(W,\LLL_W)$ is canonical.

\begin{remark}
\label{rem-ruled}
According to \cite[Cor. 2.14]{reid} all the components of $B$ are
birationally ruled surfaces.
\end{remark}

\begin{lemma}
\label{lemma-last-7}
In notation~\xref{notation_1}, we can take a modification $\phi$
so that
\begin{enumerate}
\item
the pair $(W, \HHH_W)$ is canonical, and
\item
the linear system $\HHH_W$ is nef.
\end{enumerate}
\end{lemma}

\begin{proof}
First we consider case (0). According to \cite[Th 2.11]{reid}
there exists a blowup $\phi_1\colon V_1\to V$ such that
\begin{enumerate}
\item
the variety $V_1$ normal;
\item
the proper transform $H_1$ of a general divisor $H\in \HHH$ is
Cartier and has only Du Val singularities;
\item
$K_{V_1}=\phi_1^*K_{V}$ (in particular, $K_{V_1}$ is Cartier and
$V_1$ has only canonical singularities).
\end{enumerate}
By Inversion of Adjunction \cite[3.3, 9.5]{Sh} the pair
$(V_1,\HHH_1)$ is purely log terminal. Since the linear system
$\HHH_1$ consists of Cartier divisors, the pair $(V_1,\HHH_1)$ is
canonical.

Consider a terminal $\QQ$-factorial modification $g_1\colon W\to
V_1$ of $V_1$. The composition
\[
\phi\colon W\stackrel{g_1}{\longrightarrow} V_1
\stackrel{\phi_1}{\longrightarrow} V.
\]
is a terminal $\QQ$-factorial modification for $V$. In this case,
$K_W=g_1^*K_{V_1}$ and $\HHH_W=g_1^*\HHH_1-E$, where $E\ge 0$.
Thus, $K_W+\HHH_W=g_1^*(K_{V_1}+\HHH_1)-E$. Since the pair
$(V_1,\HHH_1)$ is canonical, $E=0$ and so is $(W,\HHH_W)$. For
another terminal modification $\phi'\colon W'\to V$ the map
$W\dashrightarrow W'$ is an isomorphism in codimension $1$, so it
is crepant. By Lemma~\xref{lemma-crep-otobr-1} the pair
$(W',\HHH')$ is also canonical.

Since the linear system $\HHH_W$ has no fixed components, there
exists only a finite number of curves $C_i\subset W$ having
negative intersection numbers with $\HHH_W$ and these curves are
contained in $\Bs \HHH_W\subset \phi^{-1}(P)$. Apply the
$\HHH_W$-Minimal Model Program to $W$ over $V$. After a finite
number of $\HHH_W$-flops we get a new terminal modification
$\phi'\colon W'\to V$ with nef linear system $\HHH_{W'}$.
Replacing $W$ with $W'$ we may assume that $\HHH_{W}$ is nef.

In cases (A) and (B) we show that the linear system $\HHH_W$ has
no base points and fixed components. In particular, the pair
$(W,\HHH_W)$ is terminal.

Consider case (A). Since $\HHH$ is a linear system of hyperplane
sections passing through $\Gamma$, for the proof of the lemma it
is sufficient to show only that the morphism $\phi$ is decomposed
through the blowup of $\Gamma$ as a reduced subscheme in $V$.

Since every crepant blowup can be extended to a terminal
$\QQ$-factorial modification, our assertion is an immediate
consequence of the following.

\begin{claim}
Let $V$ be a threefold with cDV singularities and let
$\Gamma\subset \Sing(V)$ be an one-dimensional irreducible
component. Assume that the curve $\Gamma$ is smooth. Then the
blowup $\sigma\colon \tilde V\to V$ of the curve $\Gamma$ as a
reduced subscheme is crepant. In particular, $\tilde V$ is normal
and has only cDV singularities.
\end{claim}

\begin{proof}
The problem is local. Hence we may assume that $V$ is an analytic
neighborhood of some point $P\in \Gamma$. Thus we may assume that
$V$ is a hypersurface singularity, given by the equation
\[
\psi_0(x,y,z)+t\psi_1(x,y,z)+t^2\psi_2(x,y,z)+\cdots=0,
\]
where $\psi_0(x,y,z)=0$ is an equation of a Du Val singularity,
and the curve $\Gamma$ is given by the equations $x=y=z=0$. Since
$V$ is singular along $\Gamma$, $\mult_{(0,0,0)}\psi_k\ge 2$ and
since $V$ has only cDV singularities, $\mult_{(0,0,0)}\psi_0= 2$.
Then in the chart $x\neq 0$ the variety $\tilde V$ is given by the
equation
\[
\psi_0(x,yx,zx)x^{-2}+t\psi_1(x,yx,zx)x^{-2}+
t^2\psi_2(x,yx,zx)x^{-2}+\cdots=0.
\]
Here $\psi_0(x,yx,zx)x^{-2}$ is not divisible by $x$ and
$\psi_0(x,yx,zx)x^{-2}=0$ is an equation of a smooth or Du Val
point. This shows that $\tilde V$ has only cDV singularities in
our chart and the fiber over the point $P$ is one-dimensional. By
symmetry the same is true in other charts. Therefore all the
components of the exceptional divisor are surjectively mapped to
$\Gamma$. It is easy to see by the Adjunction Formula that the
morphism $\sigma$ is crepant.
\end{proof}
Case (B) is considered in a similar way. Lemma~\xref{lemma-last-7}
is proved.
\end{proof}

Now run the $K_W+\HHH_W$-Minimal Model Program. On each step
relation~\eqref{eq-HB} is kept, so the log divisor $K+\HHH\equiv
-B$ can not be nef. At the end we get a canonical pair
$(X,\HHH_X)$ and $K_X+\HHH_X$-negative contraction $f\colon X\to
Z$ to a lower-dimensional variety $Z$. Moreover, $X$ has only
$\QQ$-factorial terminal singularities (however $X$ is not
necessarily Gorenstein). By Lemma~\xref{lemma-crep-otobr-1} the
pair $(X,\LLL_X)$ is canonical and $\LLL_X\subset |-K_X|$. The
Fano threefold $V=V_{2g-2}\subset \PP^{g+1}$ is the image of $X$
under the birational map defined by the linear system $\LLL_X$.
According to~\eqref{eq-HB} we have
\[
K_X+\HHH_X+B_X\sim 0.
\]
By Lemma~\xref{lemma-predv-image} the image $\Phi_{\HHH_X}(X)$ is
three-dimensional. In particular, the linear system $\HHH_X$ is
not a pull-back of a linear system on $Z$, i.e., $\HHH_X$ is ample
over $Z$.

Everywhere below we assume that $-K_V^3\ge 72$. Then according to
\eqref{lemma-g-dim} we have $\dim |-K_V|\ge 38$ and $\dim \HHH\ge
35$. Hence
\begin{gather}
\label{eq-main-K}
\dim |-K_X|\ge \dim \LLL_X =\dim |-K_V|\ge 38,
\\
\label{eq-main-H}
\dim |H|\ge \dim \HHH\ge 35,
\end{gather}
where $H\in \HHH_X$ is a general divisor.

For $Z$ there are only the following possibilities: a) $Z$ is a
point, b) $Z$ is a curve, and c) $Z$ is a surface. In case b), a
general fiber $X_\eta$ is a nonsingular del Pezzo surface and
divisors $H|_{X_\eta}$ and $-(K_{X_\eta}+H|_{X_\eta})$ are ample.
Therefore, $X_\eta\simeq\PP^2$ or $\PP^1\times\PP^1$. Below, in
next sections, we consider cases according to the dimension of $Z$
and the type of the fiber $X_{\eta}$. In each case we obtain
contradiction with~\eqref{eq-main-K}-\eqref{eq-main-H} or show
that in~\eqref{eq-main-K} equalities hold. In the last case
$|-K_X|=\LLL_X$ and the birational map $X \dashrightarrow V$ is
given by the complete linear system $|-K_X|$. Then it is easy to
show that, for $V$, there are only two possibilities as in Example
\xref{ex-main}.

\section{Case: $Z$ is a point}
\begin{pusto}
\label{not-point}
In this section we consider the case, when $Z$ is a point. Then
$\rho(X)=1$ and $X$ is a Fano threefold with $\QQ$-factorial
terminal singularities. In this situation $\Pic X\simeq \ZZ$ (see,
e.g., \cite[Prop. 2.1.2]{IP}). Let $G$ be an ample Cartier divisor
that generates $\Pic X$. We can write $-K_X\equiv rG$ and
$\HHH_X\sim aG$ for some $r\in \QQ$, $r>0$ and $a\in \NN$. Such an
$r$ is called the \emph{Fano index} of $X$. Since the divisor
$-(K_X+\HHH_X)$ is ample, we have $r>a\ge 1$.
\end{pusto}

\begin{proposition}
Notation as in~\xref{not-point}. Then $\dim |-K_X|\le 34$.
\end{proposition}
\begin{proof}
If $X$ has only Gorenstein singularities, then by Theorem
\xref{th-Namik} we have $\dim |-K_X|\le 34$. Therefore $X$ has at
least one point of index $>1$. Below we use the following theorem.

\begin{theorem}[\cite{S}]
\label{th-Fan-pt}
Let $X$ be a Fano threefold with terminal singularities. Assume
that $X$ has at least one point of index $>1$ and the Fano index
$r>1$. Then there is one of the following embeddings of $X$ into a
weighted projective space \textup(we use numeration of
\cite{S}\textup)\textup:
\begin{enumerate}
\renewcommand\labelenumi{[\arabic{enumi}]}
\item
$X=X_6\subset \PP(1,1,2,3,i)$, $i=2,3,4,5,6$, $r=1 + 1/i$;
\item
$X=X_4\subset \PP(1,1,1,2,i)$, $i=2,3$;
\item
$X=X_3\subset \PP(1,1,1,1,2)$;
\setcounter{enumi}{4}
\item
$X=\PP(1,1,1,2)$.
\end{enumerate}
\end{theorem}

Recall that a weighted projective space $\PP=\PP(w_0,\dots,w_n)$
is said to be \emph{normalized} if
$\gcd(w_0,\dots,w_{j-1},w_{j+1},\dots,w_n)=1$ for $j=0,\dots,n$.
Under this condition, the canonical divisor is computed by the
formula $\OOO_{\PP}(K_{\PP})=\OOO_{\PP}(-\sum d_j)$. All weighted
projective spaces $\PP$ in Theorem~\xref{th-Fan-pt} are
nonsingular in codimension $2$. Therefore, in cases [1]-[3], the
standard the Adjunction Formula $K_X=(K_{\PP}+X)|_X$ holds. From
the exact sequence
\[
0\longrightarrow \OOO_{\PP}(-K_{\PP}-2X) \longrightarrow
\OOO_{\PP}(-K_{\PP}-X) \longrightarrow \OOO_{X}(-K_X)
\longrightarrow 0
\]
and the vanishing $H^1(\PP,\OOO_{\PP}(-K_{\PP}-2X))=0$ it is easy
to compute the dimension of $|-K_X|$. Here we use the fact that
$H^0(\PP,\OOO(d))$ is naturally isomorphic to the component of
degree $d$ of the polynomial ring $\CC[x_0,\dots,x_n]$ with
graduation given by $\deg x_j=w_j$. Then in [1]-[3] the maximal
value $\dim |-K_X|=30$ is achieved in case [1] for $i=6$. In case
[5] we have $-K_X\sim \OOO(5)$ and $\dim |-K_X|=33$.
\end{proof}

\section{Case: $X_{\eta}\simeq\PP^2$}
\label{sect-p2-2}
Since $H^1(Z,\OOO_Z)=H^1(X,\OOO_X)=0$ (by the Kawamata-Viehweg
Vanishing Theorem), $Z\simeq\PP^1$. Let $H\in \HHH_X$ be a general
element. The divisor $-(K_{X_\eta}+H|_{X_\eta})$ is ample, so
$H|_{X_{\eta}}\simeq\OOO_{\PP^2}(1)$ or $\OOO_{\PP^2}(2)$.

\begin{pusto}
\label{not-curve-p2-o1}
First we consider the case $H|_{X_{\eta}}\simeq\OOO_{\PP^2}(1)$.
The following fact is well-known (see, e.g., \cite{Fujita1}). For
convenience of the reader we give the proof.
\end{pusto}

\begin{lemma}
\label{lemma-p2-p2-bundle-62}
If $H|_{X_{\eta}}\simeq\OOO_{\PP^2}(1)$, then $f$ is a
$\PP^2$-bundle.
\end{lemma}
\begin{proof}
Let $S=g^*P$ be an arbitrary fiber. Then $H^2\cdot S=H^2\cdot
X_\eta=1$. Therefore $S$ is reduced and irreducible. Since the
morphism $f$ is flat, the function $\chi(X_z,\OOO_{X_z}(H))$ is
locally constant. Thus, $\chi(S,\OOO_{S}(H))=\chi
(\OOO_{\PP^2},\OOO_{\PP^2}(1))=3$. On the other hand,
$h^2(S,\OOO_S(H))=h^0(S,\omega_S \otimes \OOO_S(-H))=0$. Hence,
$h^0(S,\OOO(H))\ge 3$. Recall the definition of $\Delta$-genus of
a polarized variety $(Y,\MMM)$ (see \cite{Fujita}):
\[
\Delta(Y,\MMM)=\dim Y+ \MMM^{\dim Y}-h^0(Y,\MMM).
\]
It is known that $\Delta(Y,\MMM)\ge 0$ and in case
$\Delta(Y,\MMM)= 0$ the variety $Y$ is normal and the sheaf $\MMM$
is very ample \cite{Fujita}. In our case we have
$\Delta(S,\OOO_S(H))\le 0$ and $\OOO_S(H)^2=1$. Hence,
$S\simeq\PP^2$.
\end{proof}

\begin{proposition}
\label{prop-main-p2-bundle}
Let $X\to Z=\PP^1$ be a $\PP^2$-bundle. Assume that the pair
$(X,|-K_X|)$ is canonical. Then $\dim |-K_X|\le 38$. Moreover, if
$\dim |-K_X|= 38$, then
$X=\PP(\OOO_{\PP^1}(6)\oplus\OOO_{\PP^1}(2)\oplus\OOO_{\PP^1})$
and the anti-canonical image of $X$ is the weighted projective
space $\PP(6,4,1,1)$.
\end{proposition}

\begin{pusto}
We may assume that $X=\PP_{\PP^1}(\EEE)$, where
\[
\EEE=\OOO_{\PP^1}(d_1)\oplus\OOO_{\PP^1}(d_2)\oplus\OOO_{\PP^1}(d_3),\quad
d_1\ge d_2\ge d_3= 0.
\]
Denote the class of the fiber of $f$ by $F$ and the tautological
divisor by $M$, i.e., a divisor such that $\OOO_X(M)\simeq
\OOO_{\PP(\EEE)}(1)$. Put $d:=\sum d_i$. The relative Euler
sequence gives us
\[
-K_X\sim 3M+(2-d)F.
\]
\end{pusto}
Further, put $\EEE':=\OOO_{\PP^1}(d_1)\oplus\OOO_{\PP^1}(d_2)$.
Then
\begin{multline}
\label{eq-p2-dim-k}
H^0(X,\OOO_X(-K_X))= H^0\bigl(\PP(\EEE),\OOO_{\PP(\EEE)}(3)\otimes
f^*\OOO_{\PP^1}(2-d)\bigr)\simeq
\\
H^0(\PP^1,S^3\EEE\otimes \OOO_{\PP^1}(2-d))=
H^0(\PP^1,S^3\EEE'(2-d))\oplus
\\
H^0(\PP^1,S^2\EEE'(2-d))\oplus H^0(\PP^1,\EEE'(2-d))\oplus
H^0(\PP^1,\OOO_{\PP^1}(2-d)),
\end{multline}
where
\[
H^0(\PP^1,S^m\EEE'(2-d))=\bigoplus_{i=0}^m
H^0(\PP^1,\OOO_{\PP^1}(id_1+(m-i)d_2+2-d)).
\]
On the other hand, the pair $(X,|-K_X|)$ is canonical (because
$|-K_X|\supset \LLL_X$). In particular, the linear system $|-K_X|$
has no fixed components and its general member has at worst
isolated singularities.

It is clear that $H^0(X,M-d_1F)\simeq H^0(\PP^1,\EEE(-d_1))\neq
0$. Therefore there exists an effective divisor $D\sim M-d_1F$.
Since the divisor $M$ is nef, $-K_X\cdot D\cdot M\ge 0$. Using
Hirsch formula $M^2=dM\cdot F$ we can compute $-K_X\cdot D\cdot M$
and obtain
\[
2d_2+2-d_1\ge 0.
\]
If now $d_2=0$, then $d_1=d\le 2$ and by~\eqref{eq-p2-dim-k} we
get $\dim |-K_X|= 29$. Assume that $d_2>0$. Consider the
decomposition $\EEE=\EEE'\oplus \OOO_{\PP^1}$ and let $C\subset X$
be the section, corresponding to the surjection $\EEE\to \OOO$. A
general section $L\in |-K_X|=|3M+(2-d)F|$ has multiplicity $\le 1$
at the general point of $C$ and the curve $C$ is given by
vanishing of all sections in $H^0(\PP^1,\EEE')\subset
H^0(\PP^1,\EEE)$. Therefore, in the sum~\eqref{eq-p2-dim-k}, the
term $H^0(\PP^1,\EEE'(2-d))$ does not vanish. This means that
$d_1+2-d\ge 0$. Thus, $d_2\le 2$ and $d_1\le 2d_2+2\le 6$. Using
these inequalities we get the following possibilities for
$(d_1,d_2)$:
\[(1, 1), (2, 1), (2, 2), (3, 1), (3, 2),
(4, 1), (4, 2), (5, 2), (6, 2).
\]
By~\eqref{eq-p2-dim-k} we can immediately compute that $\dim
|-K_X|\le 38$. The equality holds only for $(d_1,d_2)=(6, 2)$. The
last statement of Proposition~\xref{prop-main-p2-bundle} follows
from the Fano construction (see \cite[Ch. 4, Remark 4.2]{Isk0}).
This proves Proposition~\xref{prop-main-p2-bundle}.

\begin{pusto}
\label{not-curve-p2-o2}
From now on to the end of this section we assume that
$X_{\eta}\simeq\PP^2$ and $H|_{X_{\eta}}\simeq\OOO_{\PP^2}(2)$.
\end{pusto}

We prove the following.
\begin{proposition}
\label{prop-p2-sing-crep}
In notation~\xref{not-curve-p2-o2} there exists a
$K_X+\LLL$-crepant birational map $X\dashrightarrow X_0$ onto a
$\PP^2$-bundle $X_0$ over $Z$.
\end{proposition}

This implies Theorem~\xref{th-main} in case
\xref{not-curve-p2-o2}. Indeed, by Proposition
\xref{prop-main-p2-bundle} we have $\dim \LLL=\dim \LLL_{X_0}\le
\dim |-K_{X_0}|\le 38$. Moreover, if equalities hold, then images
of $X$ and $X_0$ under birational maps defined by linear systems
$\LLL$ and $|-K_{X_0}|$ coincide.

\begin{lemma}
\label{lemma-p2-2-1}
Let $S=f^*z_0$, $z_0\in Z$ be a reduced fiber. Then
\begin{enumerate}
\item
$S$ is a normal surface;
\item
$S\simeq \PP^2$ or $\cone_4$;
\item
the pair $(X,S)$ is purely log terminal.
\end{enumerate}
\end{lemma}
\begin{proof}
As in the proof of Lemma~\xref{lemma-p2-p2-bundle-62} we have
$h^0(S,\OOO(H))\ge 6$ and $\Delta(S,\OOO_S(H))= 0$. Therefore the
surface $S$ is normal and there is an embedding $S\hookrightarrow
\PP^5$ whose image is a surface of degree $4$ (because
$\OOO_S(H)^2=4$). It is well-known that in this situation for $S$
there are only possibilities in (ii). The statement (iii) follows
now by Inversion of Adjunction \cite[3.3, 9.5]{Sh} (because $X$ is
nonsingular in codimension $2$).
\end{proof}

\begin{corollary}
\label{lemma-p2-2-2}
The bundle $f\colon X\to Z$ has no multiple fibers.
\footnote{There is a gap in the proof. 
See [math.AG/0604468, Lemma 5.9] for corrections.}
\end{corollary}
\begin{proof}
Assume that $f^*z_0=mS$ is a fiber of multiplicity $m>1$. Consider
$X$ as a small neighborhood of $S$ and consider the base change
\[
\xymatrix{
X\ar[d]^f&X'\ar[d]^{f'}\ar[l]_{\varphi}
\\
Z& Z'\ar[l]^{\psi}
}
\]
where $\psi$ is a covering locally given by $t\to t^m$ and $X'$ is
the normalization of the dominant component of $X\times_Z Z'$.
Then $\varphi$ is a finite cyclic covering of degree $m$ which is
\'etale outside of $\Sing X$. Put $S':=\varphi^{-1}(S)$. Then $S'$
is a Cartier divisor and coincides with the scheme-theoretical
fiber: $S'=f'^*z_0'$, where $z_0'=\psi^{-1}(z_0)$. Since $\Sing X$
is a finite set, components $S'$ can intersect each other only in
a finite number of points $\varphi^{-1}(\Sing X)$. This
contradicts Lemma~\xref{lemma-vspom-term-CM} below. Therefore the
fiber $S'$ is irreducible and reduced. By Lemma
\xref{lemma-p2-2-1}, $S'\simeq\PP^2$ or $\cone_4$. However a
cyclic group acting on $\PP^2$ or $\cone_4$ has a curve of fixed
points. A contradiction shows that $m=1$.
\end{proof}

\begin{lemma}
\label{lemma-vspom-term-CM}
Let $P\in Y$ be a three-dimensional terminal singularity and let
$D'$, $D''$ be effective Weil divisors such that $D=D'+D''$ is
Cartier. Then $\Supp D'\cap \Supp D''\ge 1$.
\end{lemma}
\begin{proof}
Assume that $\Supp D'\cap \Supp D''=\{P\}$. It is clear that we
can replace $P\in Y$ with its canonical covering. Thus we may
assume that $P\in Y$ is an isolated hypersurface singularity. Then
$D$ is a locally complete intersection and $\Supp D'\cap \Supp
D''=\{P\}$ which is impossible.
\end{proof}

\begin{corollary}
If $S\simeq \PP^2$, then $X$ is nonsingular along $S$. If
$S\simeq\cone_4$, then $X$ has exactly one singular point which is
the vertex of the cone $\cone_4$ and this point is analytically
isomorphic to
\begin{equation}
\label{eq-p2--2-sing-point}
\{x_1x_2+x_3^2+x_4^n=0\}/\muu_2(1,1,1,0),\quad n\ge 1.
\end{equation}
\end{corollary}
\begin{proof}
If the variety $X$ is singular at some point $P\in S$, then so is
the surface $S$ (because $S$ is a Cartier divisor). Assume that
$S\simeq\cone_4$. Then the (Gorenstein) index of the vertex $O\in
S$ is equal to the index of the point $O\in X$. Therefore $O\in X$
is a point of index $2$ (in particular, this point is singular).
Now we consider $X\ni O$ as a small neighborhood and consider the
canonical $\muu_2$-covering $X'\to X$ near $O$. This induces the
covering $S'\to S$, where $S'\simeq \CC^2/\muu_2(1,1,1)$, i.e.,
$S'$ is a singularity of type $A_1$. Since $S'\subset X'$ is a
Cartier divisor, we may assume that $X'\ni O'$ is given by the
equation $x_1x_2+x_3^2+x_4\phi(x_1,x_2,x_3,x_4)$ and $\muu_2$ acts
on $x_1,x_2,x_3$ by $x_i\to -x_i$. The rest follows by the
classification of terminal singularities.
\end{proof}

\begin{lemma}
\label{lemma-p2-singpot-blowup}
Let $X\ni O$ be a singular point of the form
\eqref{eq-p2--2-sing-point}, let $\sigma\colon \tilde X\to X$ be
the weighted blowup with weights $\frac12(1,1,1,2)$ and let $E$ be
the exceptional divisor. Then
\begin{enumerate}
\item
$E\simeq \cone_4$ for $n\ge 2$ and $E\simeq \PP^2$ for $n=1$;
\item
$a(E)=1/2$;
\item
if $n=1$, then the variety $\tilde X$ is nonsingular, if $n\ge 2$,
then $\tilde X$ has exactly one singular point at the vertex
$O_4\in \cone_4$ which \textup(up to analytic isomorphism\textup)
has the form
\begin{equation}
\label{eq-p2--2-sing-point-11}
\{x_1x_2+x_3^2+x_4^{n-1}=0\}/\muu_2(1,1,1,0).
\end{equation}
\end{enumerate}
\end{lemma}
\begin{proof}
The divisor $E$ is given in $\PP(1,1,1,2)$ by the equation
$x_1x_2+x_3^2=0$ for $n\ge 2$ and $x_1x_2+x_3^2+x_4^2=0$ for
$n=1$. This implies (i). Computations of discrepancies in (ii) is
the standard toric technique. We prove (iii). The variety $\tilde
X$ is covered by four affine charts $U_i$. The map
$X\dashrightarrow U_1\simeq \CC^4$ is given by the following
formulas
\[
(x_1,x_2,x_3,x_4)\longrightarrow
(x_1^{1/2},x_2x_1^{1/2},x_3x_1^{1/2},x_4x_1)
\]
In this chart $\tilde X=\{x_2+x_3^2+x_4^nx_1^{n-1}\}$ is
nonsingular. Computations for $U_2$ and $U_3$ are completely the
same.

The map $X\dashrightarrow U_4\simeq \CC^4/\muu_2(1,1,1,0)$ has the
form
\[
(x_1,x_2,x_3,x_4)\longrightarrow
(x_1x_4^{1/2},x_2x_4^{1/2},x_3x_4^{1/2},x_4)
\]
Thus, in this chart, $\tilde
X=\{x_1x_2+x_3^2+x_4^{n-1}\}/\muu_2(1,1,1,0)$ has exactly one
singular point at the origin. The rest is obvious.
\end{proof}

To finish the proof of Proposition~\xref{prop-p2-sing-crep} we
consider the weighted blowup $\sigma\colon \tilde X\to X$ as in
Lemma~\xref{lemma-p2-singpot-blowup}. Then
\[
K_{\tilde X}+\tilde S=\sigma^*(K_X+S)-\frac12E, \quad
\sigma^*S=\tilde S+E
\]
(because $S$ is Cartier and the pair $(X,S)$ is purely log
terminal). Therefore the pair $(\tilde X,\tilde S+\frac12E)$ is
also purely log terminal. Hence $\tilde S$ is a normal surface.
Note that the divisor
\[
K_{\tilde S}=(K_{\tilde X}+\tilde S)|_{\tilde S}=-\frac12
E|_{\tilde S}+\sigma^*(K_X+S)|_{\tilde S}
\]
is $\sigma|_{\tilde S}$-ample. Therefore $\tilde S$ is the minimal
resolution of the singularity $S\simeq\cone_4$, so $\tilde
S\simeq\FF_4$. Since the surface $\tilde S$ is nonsingular, so is
the variety $\tilde X$ along $\tilde S$. It is clear that $\tilde
S\cap E$ is the minimal section $\Sigma$ on $\tilde S\simeq
\FF_4$. Let $l$ be a fiber of the projection $\tilde S\simeq
\FF_4\to \PP^1$. Note that the intersection $\tilde S\cap E$ is
reduced at the general point (because the pair $(\tilde X,\tilde
S+\frac12E)$ purely log terminal). Then $\tilde S\cdot
l=\sigma^*S\cdot l-E\cdot l=-\Sigma\cdot l=-1$. According to the
contractibility criterion there exists a birational contraction
$\varphi\colon \tilde X\to X_{n-1}$ over $Z$ which contracts the
surface $\tilde S$. The curve $\varphi(\tilde S)$ is contained in
the nonsingular locus of $X_{n-1}$. The variety $X_{n-1}$ has
exactly one singular point and this point has the form
\eqref{eq-p2--2-sing-point-11}. Continuing the process we get the
following series of birational transformations over $Z$:
\[
X=X_n\dashrightarrow X_{n-1}\dashrightarrow \cdots \dashrightarrow
X_1 \dashrightarrow X_0
\]
On the last step the variety $X_0$ is nonsingular. Thus, $X_0\to
Z$ is a $\PP^2$-bundle. This proves Proposition
\xref{prop-p2-sing-crep}.

\section{Case: $X_{\eta}\simeq\PP^1\times\PP^1$}
\begin{pusto}
In this section we assume that $f\colon X\to Z$ is an extremal
Mori contraction with general fiber $X_{\eta} \simeq \PP^1
\times\PP^1$. Let $H\in \HHH_X$ be a general element.
\end{pusto}

\begin{lemma}[cf. {\cite[3.5]{Mori}}]
There exists an embedding of $X$ into a $\PP^3$-bundle over $Z$
such that every fiber $X_{\eta}$ is a reduced irreducible quadric
\textup(in particular, $X$ is Gorenstein and every Weil divisor on
$X$ is Cartier\textup).\footnote{There is a gap in the proof. 
See [math.AG/0604468, Lemma 5.9] for corrections.}
\end{lemma}
\begin{proof}[Sketch of the proof]
Let $f^*z_0$ be an arbitrary fiber. Similar to Lemma
\xref{lemma-p2-2-1} and Corollary~\xref{lemma-p2-2-2} one can
prove that $f^*z_0$ is an irreducible reduced normal surface
isomorphic to $\PP^1\times\PP^1$ or a quadratic cone $\cone_2$.
Therefore the divisor $H$ is relatively very ample over $Z$ and
defines the desired embedding.
\end{proof}

\begin{pusto}
\label{notat-Q-2}
Thus there is an embedding $X\hookrightarrow \PP$ over $Z$, where
\[
\PP=\PP(\EEE),\qquad \EEE=\bigoplus_{i=1}^{4} \OOO_{\PP^1}(d_i)
\]
and $X_{\eta}\subset \PP_{\eta}$ is a quadric. Put $d=\sum d_i$.
We can chose $d_i$ so that $d_1\geq d_2\geq d_3\geq d_4=0$. Let
$M$ be the tautological divisor on $\PP$ and let $F$ be a fiber of
the projection $\pi\colon \PP\to Z$.
\end{pusto}
We prove the following.
\begin{proposition}
\label{prop-p1-p1-m}
If, in notation~\xref{notat-Q-2}, $\dim |-K_X|\ge 38$, then there
exists a $K_X+\LLL$-crepant birational map of $X$ onto a
$\PP^2$-bundle over $Z$.
\end{proposition}
As in \S~\xref{sect-p2-2} this is sufficient for the proof of
Theorem~\xref{th-main} in the case $X_\eta\simeq
\PP^1\times\PP^1$. Recall well-known facts which are consequences
of the Hirsch formula and relative Euler sequence.

\begin{lemma}
\label{lemma-Pic-P-int-bas}
\begin{enumerate}
\item
$\Pic \PP=\ZZ\cdot M\oplus \ZZ\cdot F$;
\item
$M^4=d$,\qquad $M^3\cdot F=1$,\qquad $F^2\equiv 0$;
\item
$-K_{\PP}=4M+(2-d)F$.
\end{enumerate}
\end{lemma}

Thus $X\sim 2M+r F$ for some $r\in \ZZ$. Put $G:=M|_X$ and
$Q:=F|_X$. By the Adjunction Formula
\[
-K_X=2G+(2-d-r)Q.
\]
\begin{corollary}
\label{cor-Pic-Q}
In the above notation, we have $\Pic X\simeq \ZZ\cdot
G\oplus\ZZ\cdot Q$.
\end{corollary}

\begin{lemma}
We may assume that $d+r\ge 3$.
\end{lemma}
\begin{proof}
If $d+r<2$, then the divisor $-K_X$ is ample (because $|G|$ and
$|Q|$ are base point free linear systems). By the classification
of three-dimensional nonsingular Fano threefolds we have $\dim
|-K_V|\le \dim |-K_X|\le 34$.

Now let $d+r=2$. Then $-K_X=2G$. In this situation a general
member of the linear system $|G|$ has only Du Val singularities
(see, e.g., \cite{Shin}). By the Adjunction Formula $-K_G$ is nef
and big. Since $\dim |-K_X|\ge 38$, as in~\eqref{lemma-g-dim} we
have $-K_X^3\ge 72$. Hence, $K_G^2=-\frac18K_X^3\ge 9$. On the
other hand, applying Noether formula to the minimal resolution
$\tilde G\to G$ we get $K_G^2\le 9$. Moreover, equality holds only
if $G=\tilde G\simeq \PP^2$ (cf.~\xref{intro-del-Pezzo}). But the
last contradicts the fact that $G\simeq \PP^2$ has no nontrivial
morphisms to a curve.
\end{proof}

We can write $-K_X\sim H+B$, where $B$ is an effective divisor. By
Corollary~\xref{cor-Pic-Q}, $H\sim G+\alpha Q$ and $B\sim
G-(d+r+\alpha-2)Q$ for some $\alpha\in \ZZ$.

\begin{lemma}
\label{lemma-ocenka-6}
$d+2r\le 6$.
\end{lemma}
\begin{proof}
Since the divisor $H$ is nef,
\[
0\le (-K_X)\cdot B\cdot H=2(6-d-2r).
\]
This proves the statement.
\end{proof}

\begin{lemma}
$r<0$.
\end{lemma}
\begin{proof}
Assume that $r\ge 0$. Then $d\le 6$. Since $R^i
f_*(\OOO_{X}(H))=0$ for $i>0$ (see Lemma~\xref{lemma-Kaw-Vieh}),
\[
h^0(\OOO_{X}(H))=h^0(\OOO_{\PP^1}(\EEE(\alpha))).
\]
If $\alpha\le 0$, then $h^0(\OOO_{X}(H))\le
h^0(\OOO_{\PP^1}(\EEE))=d+4\le 10$. If $\alpha\ge 0$, then we have
$\alpha \le 2+d_1-d-r$ (because $B$ effective) and
\[
h^0(\OOO_{X}(H))= d+4+4\alpha\le 12 +4d_1 -3d- 4r\le 24.
\]
In both cases we have a contradiction with our assumption
\eqref{eq-main-H}.
\end{proof}

The linear system $|M|$ has no fixed components and base points
and defines a birational contraction of $\PP$ onto a cone in
$\PP^{d+3}$. The subvariety $C:=\PP(\oplus_{d_i=0}
\OOO_{\PP^1}(d_i))$ is contracted to the vertex of the cone. This
subvariety is swept out by curves $\Gamma$ such that $M\cdot
\Gamma=0$. Since $r<0$, $C$ is contained in $X$. If $\dim C=2$,
then we have a contradiction with $\rho(X/Z)=1$. Therefore $C$ is
a curve.

\begin{proof}[Proof of Proposition~\xref{prop-p1-p1-m}]
Since $d+r>2$, we have $-K_X\cdot C<0$. Therefore, $C\subset \Bs
|-K_X|$. Let $\sigma\colon (\tilde \PP\supset \tilde X)\to
(\PP\supset X)$ be the blowup of $C$. Then $-K_{\tilde X}$ is an
ample over $X$ Cartier divisor and the exceptional divisor
$E\subset \tilde X$ of the contraction $\sigma\colon \tilde X\to
X$ is irreducible. Since the variety $\tilde X$ is nonsingular in
codimension $1$, it is normal. It is clear that $a(E,\HHH)=0$.
Therefore the contraction $\sigma$ is crepant: $K_{\tilde
X}+\tilde\HHH=\sigma^*(K_X+\HHH)$. Thus $X'$ has only
$\QQ$-factorial canonical singularities and $\rho(X'/Z)=2$. Note
that the divisor $\sigma^*G-E$ is nef over $Z$ (i.e., it is nef on
the fibers of $f\circ\sigma\colon \tilde X\to Z$). A general fiber
$\tilde Q$ of the morphism $f\circ\sigma$ is isomorphic to a
blowup of a point on $Q\simeq\PP^1\times\PP^1$. Therefore on
$\tilde Q$ there are exactly two curves $\tilde C_1$, $\tilde C_2$
having intersection number $0$ with $\sigma^*G-E$. These curves
are is proper transforms of generators $C_1$, $C_2$ passing
through the point $Q\cap C$. Therefore the divisor $\sigma^*G-E$
defines a supporting function for an extremal ray on $\Mori(\tilde
X/Z)$ and the linear system $|n(\sigma^*G-E)|$ for $n\gg 0$ gives
an extremal birational contraction $\psi\colon \tilde X\to X'$
over $Z$. The morphism $\psi$ contracts the divisor swept out by
$\tilde C_i$ to a curve. This shows that a general fiber of the
morphism $X'\to Z$ is isomorphic to $\PP^2$. Further, we can write
$\sigma^*G-E=\psi^*G'$ for some ample over $Z$ Cartier divisor
$G'$ on $X'$. Let $Q'$ be an arbitrary fiber of the morphism
$X'\to Z$. Then $(G')^2\cdot Q'=(\sigma^*G-E)^2\cdot \tilde G=1$.
As in \S~\xref{sect-p2-2} one can prove that $Q'\simeq\PP^2$.
Therefore $X'\to Z$ is a $\PP^2$-bundle. Since $C\subset \Bs
|-K_X|$, the map $X\dashrightarrow X'$ is crepant.
\end{proof}

\section{Case: $\dim Z=2$}
In this section we assume that $Z$ is surface. The proof of the
following Lemma is a particular case the proof of Theorem 11.8 in
\cite{Fujita1}. Let $H$ be a general member of the linear system
$\HHH_X$.
\begin{lemma}
The surface $Z$ is nonsingular and the contraction $f\colon X\to
Z$ is a $\PP^1$-bundle.
\end{lemma}

\begin{proof}
If $C_{\gen}$ is a general fiber, then $C_{\gen}\simeq\PP^1$,
$-K_X\cdot C_{\gen}=2$ and $H\cdot C_{\gen}=1$. Let
$C=f^{-1}(P)_{\red}$ be an arbitrary fiber with the reduced
structure and let $m:=H\cdot C$. For an ample divisor $A$ on $Z$
and $a,b\gg 0$, we consider a general divisor $S\in |aH+bf^*A|$.
By Bertini's theorem $S$ is a nonsingular surface meeting $C$
transversally at $ma$ points. Therefore the restriction
$f|_S\colon S\to Z$ is a finite morphism of degree $ma\ge a$ near
$C$. On the other hand, the degree of $f|_S$ is exactly $a$
(because $S$ meets a general fiber at $a$ points). Thus the
morphism $f|_S$ is \'etale near $C$ and the point $P\in Z$ is
nonsingular. From this we immediately infer that $f$ is flat (see
\cite[23.1]{Mat}) and every fiber $f^{-1}(P)$ is a reduced
irreducible curve arithmetic genus $0$. Hence,
$f^{-1}(P)\simeq\PP^1$ (in particular, $f$ is smooth).
\end{proof}

The surface $Z$ is rational (because, for example,
$H^1(Z,\OOO_Z)=H^1(X,\OOO_X)=0$ and $Z$ is dominated by a
component of the divisor $B$ which is a birationally ruled
surface, see Remark~\xref{rem-ruled}). Put $\EEE:=f_*\OOO_X(H)$.
Then $\EEE$ is a rank $2$ vector bundle and $X\simeq \PP(\EEE)$.

We fix the following notation and conventions:
\begin{setup}
\label{not-surf-main}
Let $Z$ be a nonsingular rational surface and let $X=\PP_Z(\EEE)$ be the
projectivization of a rank-$2$ vector bundle $\EEE$ over $Z$ such
that the pair $(X,\LLL\subset |-K_X|)$ is canonical, and let $f:X\to Z$
be the projection. We assume that the following conditions hold.
\begin{enumerate}
\item
\label{new-ass-0}
$(X,\Theta)$ is a generating $0$-pair for some boundary $\Theta$, see Definition~\ref{def:FT} (in other words, $X$ is an \textit{FT-variety} \cite{PSh});

\item
$-K_X\sim H+B$, where $H$ and $B$ are effective $f$-ample divisors and $f_*\OOO_X(H)=\EEE$;
\item
the divisor $H$ is nef and the image of the 
map the map $\Phi_{|H|}$ given by the linear
system $|H|$ is three-dimensional;

\item
the pair $(X,|-K_X|)$ is canonical.
 \end{enumerate}
\end{setup}

\begin{proposition}
\label{prop-surf-main}
In notation~\xref{not-surf-main} we have either $\dim |H|\le 34$
or $\dim |-K_X|=38$ and 
the
anti-canonical image of $X$ is the weighted projective space
$\PP(3,1,1,1)$.
\end{proposition}
Theorem~\ref{th-main} is a consequence of this fact (see~\eqref{eq-main-H}).
The rest of this section is devoted to the proof of Proposition~\ref{prop-surf-main}.

\begin{setup}
\label{not-surf-main-new}
Additionally to~\ref{not-surf-main} we assume that $\dim |H|\ge 35$.
Note that 
\[
H^i(Z,\EEE)=H^i(X,\OOO_X(H))=0\quad \text{for}\quad i>0
\]
by the Kawamata-Viehweg vanishing and~\ref{not-surf-main}\ref{new-ass-0}.
Thus 
\begin{equation}
\label{eq:new:h0E}
 h^0(Z,\EEE)=\chi(Z,\EEE)\ge 36.
\end{equation} 
\end{setup}
Suppose that the surface $Z$ contains a $(-1)$-curve. We use
notation of Lemma~\xref{lemma-surf-(-1)-curve}. By this lemma we
can construct a sequence of birational transformations
\[
\xymatrix{
X=X_1\ar[d]^{f=f_1}\ar@{-->}[r]&X_2=X'\ar[d]^{f_2}\ar@{-->}[r]&\cdots\ar@{-->}[r]&X_N\ar[d]^{f_N}
\\
Z=Z_1\ar[r]& Z_2=Z'\ar[r]&\cdots\ar[r]&Z_N
}
\]
where each square is one of transformations (iii)-(v). In case (v)
we put $Z_i=Z_{i+1}$. If on some step the invariant $n\ge 2$, then
we apply the transformation (v) decreasing $n$. If $n=0$ or $1$,
then we apply the transformation (iii) or (iv) respectively. In
these cases the $(-1)$-curve on the base is contracted and our
sequence terminates (i.e., this is possible only on the last
step).

\begin{lemma}
\label{lemma-surf-(-1)-curve-2}
In notation above for the variety $X_N/Z_N$ all conditions of
\xref{not-surf-main} hold. Moreover, $\dim |H_N|\ge \dim |H|$.
\end{lemma}

Note however that the canonical property of $(X,|H|)$ is not
preserved under our transformations.
\begin{proof}
For $n\le 1$ the map $X'\dashrightarrow X$ does not contract any
divisors. By Lemma~\xref{lemma-crep-otobr-1} the pair $(X',\LLL')$
is canonical. For $n\ge 2$, we have $K_X\cdot \Sigma<0$.
Therefore, $\Sigma\subset \Bs \LLL$ and the morphism $\sigma$ is
crepant. By Lemma~\xref{lemma-crep-otobr-2} the pair $(X',\LLL')$
is again canonical. So the canonical property of the pair $(X,\LLL)$ is
preserved.

Let $H_N$ and $B_N$ be proper transforms of $H$ and $B$ on $X_N$.
Note that the map $\psi\colon X_1 \dashrightarrow X_N$ is an
isomorphism on $X_1\setminus D$ and the image $\psi(D)$ is contained in a
fiber $\Gamma$ of the projection $f_N\colon X_N\to Z_N$. This
shows that $-K_{X_N}\sim H_N+B_N$. Consider a ``Hironaka hut''
\[
\xymatrix{
&U\ar[ld]_{\psi_1}\ar[rd]^{\psi_N}&
\\
X_1\ar@{-->}[rr]&&X_N.
}
\]
and write
\[
H_U=\psi_1^*H -\sum a_iE_i=\psi_N^*H_N-\sum b_iE_i- dD_U,
\]
where $H_U$ and $D_U$ are proper transforms of $H_1$ and $D$ on
$U$, $E_i$ are prime divisors which are exceptional for $\psi_1$,
as well as for $\psi_N$, and $a_i,\, b_i,\, d\ge 0$. Then the
divisor
\[
\psi_1^*H-\psi_N^*H_N= \sum (a_i-b_i)E_i- dD_U
\]
is nef over $X_N$. By \cite[1.1]{Sh} we have $a_i\le b_i$ and
$d\ge 0$. Let $C$ be an irreducible curve on $X_N$ and let $\tilde
C\subset U$ be an irreducible curve dominating $C$. If $C=\Gamma$,
then $H_N\cdot C>0$ (because $\rho(X_N/Z_N)=1$). Assume that
$C\neq \Gamma$. Since $\psi_N(E_i),\, \psi_N(D_U)\subset \Gamma$,
we have $E_i\cdot \tilde C\ge 0$ and $D_U\cdot \tilde C\ge 0$.
Hence
\[
\psi_N^*H_N\cdot \tilde C= \psi_1^*H\cdot \tilde C+ \sum
(b_i-a_i)E_i\cdot \tilde C+ dD_U\cdot \tilde C\ge 0.
\]
This implies immediately the numerical effectiveness of $H_N$. The
lemma is proved.
\end{proof}

\begin{corollary}
In the assumptions of~\xref{not-surf-main} and~\xref{not-surf-main-new} we may assume that $Z\simeq
\PP^2$ or $\FF_e$ with $e\neq 1$.
\end{corollary}
\begin{proof}
If $Z$ contains a $(-1)$-curve $C$, then by Lemmas
\xref{lemma-surf-(-1)-curve} and~\xref{lemma-surf-(-1)-curve-2} we
can contract it. All the properties~\xref{not-surf-main} are
preserved.
\end{proof}

\begin{lemma}
\label{lemma-p1-p2-nachalo}
Let $\Gamma\subset Z$ be a nonsingular rational curve such that
$\dim |\Gamma|>0$ and let $\EEE|_\Gamma\simeq
\OOO_{\PP^1}(d_1)+\OOO_{\PP^1}(d_2)$. Then $|d_1-d_2|\le
2+\Gamma^2$.
\end{lemma}
\begin{proof}
Let $m=|d_1-d_2|$ and $G:=g^{-1}(\Gamma)$. Then $G\simeq \FF_m$.
We have
\[
-2+m=-2-\Sigma^2=K_G\cdot \Sigma=K_X\cdot \Sigma+ G\cdot
\Sigma=K_X\cdot \Sigma+\Gamma^2.
\]
Since the linear system $|-K_X|$ has no fixed components,
$K_X\cdot \Sigma\le 0$ and $m\le 2+\Gamma^2$.
\end{proof}

From now on we assume that $\dim |H|\ge 36$, i.e.,
$h^0(Z,\EEE)=\chi(\EEE,Z)\ge 37$ (see Lemma
\xref{lemma-Kaw-Vieh}). For simplicity we put $c_i:=c_i(\EEE)$,
$i=1,2$. Then
\[
-K_X= 2 H +f^*(-K_Z-c_1),
\]
\[
H^2=H\cdot f^*c_1-f^*c_2,\qquad H^3=c_1^2-c_2.
\]
Recall the Riemann-Roch formula for rank $2$ vector bundles over a
(rational) surface $Z$:
\begin{equation}
\label{eq-surf-RR}
\chi(\EEE)=\frac12(c_1^2-2c_2-K_Z\cdot c_1)+2.
\end{equation}
Since $L\cap B$ is an effective $1$-cycle, for any nef divisor $N$
on $Z$ we have
\begin{multline}
\label{eq-surf-ocenka-1-new}
0\le -K_X\cdot B\cdot f^*N=
\\
(2H+f^*(-K_Z-c_1))\cdot (H+f^*(-K_Z-c_1))\cdot f^*N=
\\
2H^2\cdot f^*N+ 3 (-K_Z-c_1)\cdot N=
\\
2 c_1\cdot N+3 (-K_Z-c_1)\cdot N= -3K_Z\cdot N-c_1\cdot N.
\end{multline}
Hence,
\begin{equation}
\label{eq-surf-ocenka-1}
c_1\cdot N\le -3K_Z\cdot N.
\end{equation}

\subsection*{Case: $Z\simeq\PP^2$}
Taking into account natural isomorphisms $H^{2i}(\PP^2,\ZZ)\simeq
\ZZ$ we assume that $c_1$ and $c_2$ are integers.
Since $\EEE$ is nef, $c_1\ge 0$.
We claim that
$c_1\le 8$. Indeed, if $c_1=9$, then in~\eqref{eq-surf-ocenka-1}
and~\eqref{eq-surf-ocenka-1-new} equalities hold. Thus, $-K_X\cdot
B\cdot f^*N=0$ for any ample divisor $N$ on $\PP^2$. Then for a
general divisor $L\in |-K_X|$, the intersection $L\cap B$ is
composed of fibers of $f$. Therefore the divisor $-K_X$ is nef
(otherwise there exists a horizontal curve $R$ such that $L\cdot
R<0$ and $B\cdot R=(L-H)\cdot R<0$). Let $B_0$ be a horizontal
component of $B$. It is easy to see that $-K_X\cdot B_0\cdot f^*N=
-K_X\cdot (B-B_0)\cdot f^*N=0$. Hence $B=B_0$ and $B$ is
contracted by the morphism $\Phi_{|-nK_X|}$ to a point. In this
case, the vector bundle $\EEE$ is decomposable. The contradiction
shows that $c_1\le 8$.

By
\eqref{eq-surf-RR} and~\eqref{eq:new:h0E} we have
\[
c_1^2+3c_1-2c_2\ge 68.
\]
\begin{corollary}
If $\EEE$ is decomposable, then $-K_X$ is nef.
\end{corollary}

\begin{proof}
Let $\EEE\simeq \OOO(d)\oplus\OOO(d+m)$, where $m\ge 0$. Then
$d\ge 0$ (because $\EEE$ is nef). Since $c_1=2d+m$ and
$c_2=d^2+dm\ge 0$,
\[
2 d^2+2 d m+m^2+6 d+3 m \ge 68.
\]
Assume that $-K_X$ is not nef. Then $0\le 2d+m\le 8$.
Taking into account that $m\le 3$ (Lemma~\xref{lemma-p1-p2-nachalo}) we get 
a contradiction.
\end{proof}

Thus $\Phi_{\LLL_X}(X)=\PP(3,1,1,1)$ if $\EEE$ decomposable (see Proposition~\xref{prop-conic}).
Further, we assume that $\EEE$ indecomposable. 

\begin{pusto}
\label{p-8-13}
There are two possibilities:

\textbf{a) $c_1$ odd.} Put $c_1=2m-3$. Then $2\le m\le 5$ and
$2m^2-3m\ge 34+c_2$. It is easy to see that
\[
c_1(\EEE(-m))=-3,\qquad c_2(\EEE(-m))=c_2-m^2+3m\le m^2- 34<0.
\]

\textbf{b) $c_1$ even.} Put $c_1=2m-2$. Then $1\le m\le 5$ and
$2m^2-m\ge 35+c_2$. Here
\[
c_1(\EEE(-m))=-2,\qquad c_2(\EEE(-m))=c_2-m^2+2m\le m^2+m-35<0.
\]

In both cases by the Riemann-Roch formula and Serre Duality we
obtain
\[
h^0(\EEE(-m))+h^0\Bigl(\EEE(-m)\otimes\det\EEE(-m)^*\otimes
\OOO(-3)\Bigr)\ge \chi(\EEE(-m))\ge 1.
\]
Therefore, $H^0(\EEE(-m))\neq 0$. Let $s\in H^0(\EEE(-m))$ be a
nonzero section. If $s$ does not vanish anywhere, then there is an
embedding $\OOO\hookrightarrow \EEE(-m)$ and $\EEE(-m)$ is
decomposable. Let $\emptyset\neq Y\subset \PP^2$ be the zero locus
of $s$. Since $c_2(\EEE(-m))< 0$, $\dim Y=1$. Choose a general
line $\Gamma\subset\PP^2$ and let the intersection $\Gamma\cap Y$
consists of $k$ points. Put $r:=-c_1(\EEE(-m))$ ($r= 2$ or $3$).
Then $\EEE(-m)|_\Gamma=\OOO(k)\oplus\OOO(-r-k)$. By Lemma
\xref{lemma-p1-p2-nachalo} we have $2k+r\le 3$. Hence, $k=0$ and
$Y=\emptyset$, a contradiction.
\end{pusto}

\subsection*{Case: $Z\simeq\FF_e$}
\label{new-ass}
It remains to consider the following situation:
\par\medskip\noindent
$V$ is birational to a variety $X$ that is represented in a form
the projectivization $X=\PP_Z(\EEE)$ of a rank-2 vector bundle $\EEE$ on 
$Z:=\FF_e$, where $e\ge 0$ and $e\neq 1$. Let $f:X\to Z$ be the projection. 
Furthermore, may assume that the following conditions are satisfied:
\begin{enumerate}
\item
\label{new-ass-0}
$X$ is an \textit{FT-variety} \cite{PSh};

\item
\label{new-ass-1}
$-K_X\sim H+B$, where $H$ and $B$ are effective $f$-ample divisors and $f_*\OOO_X(H)=\EEE$;

\item
\label{new-ass-2}
the divisor $H$ is nef and big, and $\dim |H|\ge \dim |-K_V|-3=g-2 \ge 35$;

\item
\label{new-ass-3}
the pair $(X,|-K_X|)$ is canonical.
\end{enumerate}

We show that the situation described above does not occur.
Recall the notation.
By 
$\Sigma$ and $l$ we denote the minimal section and a fiber of $Z=\FF_e$, respectively.
Write
\[
\mathrm{c}_1(\EEE)=a\Sigma+b l,\qquad \mathrm{c}_2(\EEE)=c.
\]
where $a,b,c$ are non-negative integers such that $a\ge 0$ and $6+3e\ge b\ge ea$. Put
\[
p:=\left \lfloor \frac a 2\right\rfloor+1,\qquad q:= \left \lfloor \frac b 2\right\rfloor+1,
\qquad 
a':=a-2p,\qquad b':=b-2q.
\]
Then $-2\le a',\, b'\le -1$ and for the twisted vector bundle $\EEE':=\EEE\otimes \OOO_{\FF_e}(-p\Sigma-ql)$ we have:
\begin{equation}
\label{eq:c1c2}
\mathrm{c}_1(\EEE')=a'\Sigma+b' l,\qquad c':=\mathrm{c}_2(\EEE')=c+eap-aq-bp-ep^2+2pq.
\end{equation} 

Note that $H^i(Z,\EEE)=H^i(X,\OOO_X(H))=0$ for $i>0$ by the Kawamata-Viehweg vanishing and~\ref{new-ass-0}.
Then by the Riemann-Roch formula we have
\begin{eqnarray}
\label{eq-dim-Fn:new}
\chi(\EEE)&=& -\frac12ea(a+1)+ab+a+b-c+2= \dim |H|+1,
\\
\label{eq-dim-Fn:new2}
\chi(\EEE')&=& -\frac12ea'(a'+1)+a'b'+a'+b'-c'+2.
\end{eqnarray}
Then \eqref{eq-dim-Fn:new} and our assumption \ref{new-ass-2} imply
\begin{equation}
\label{eq-dim-Fn:new-a}
c\le -\frac12ea(a+1)+ab+a+b+3-g\le -\frac12ea(a+1)+ab+a+b-34
\end{equation} 

\begin{lemma}
\label{new:lemma:e}
$e\le 3$.
\end{lemma}

\begin{proof}
Assume that $e\ge 4$. 
Let $L\in |-K_X|$ be a general member.
Since the pair $(X,|-K_X|)$ is canonical, the surface $L$ has only Du Val singularities.
Then it is easy to see that $L$ is a Du Val K3 surface. The restriction map $f_L: L\to Z$ is generically finite of degree $2$.
Let 
\[
f|_{L}: L\overset{\varsigma}\longrightarrow \bar L\overset{\pi}\longrightarrow Z=\FF_e
\] 
be its Stein
factorization. Here the morphism $\pi$ is finite of degree $2$ and 
$\varsigma$ is a crepant birational contraction. This implies that 
$\bar L$ is a Du Val K3 surface.
By the Hurwitz formula, 
\[
\textstyle
K_{\bar L}=\pi^*\left(K_Z+\frac12 R\right), 
\]
where
$R\in |-2K_Z|$ is the branch divisor. 
Since $\bar L$ is connected and normal, $R$ is reduced.
This immediately implies that $e=4$ and $R=\Sigma+R'$, where $R'\cap \Sigma=\varnothing$.
Let $\bar\Upsilon:=\pi^{-1}(\Sigma)_{\mathrm{red}}$ and $\Upsilon:=\varsigma^{-1}(\bar\Upsilon)$.
Then the surface $\bar L$ is smooth over $\Sigma$ and the morphism $\pi$ induces an isomorphism 
$\bar\Upsilon\simeq \Sigma$.
This implies that the morphism $\varsigma$ is an isomorphism over $\Sigma$ and
$L$ is smooth over $\Sigma$ as well.

Denote $F:=f^{-1}(\Sigma)$. Then $F$ is a ruled surface $\FF_{e'}$.
By the above arguments the intersection $F\cap L$ does not contain any fibers of $f$,
hence $F\cap L$ is an irreducible. Moreover, $F\cap L=2\Upsilon$, where $\Upsilon:=\varsigma^{-1}(\bar\Upsilon)$
is a section of $F/\Sigma$.
By the adjunction formula
\[
2\Upsilon\sim -K_X|_F\sim -K_F+F|_F \sim 2\Sigma'+ (2+e'-e)l',
\]
where $\Sigma'$ and $l'$ be the minimal section and fiber of $F=\FF_{e'}$, respectively.
Since $\Upsilon$ is irreducible, we have $2+e'-e\ge 2 e'$ and so $e+e'\le 2$, a contradiction.
\end{proof}

\begin{lemma}
\label{new:lemma:c2}
We have $c'=\mathrm{c}_2(\EEE')\le -2$. Moreover, the equality holds only if
\begin{equation}
\label{new:lemma:c2:case1}
e=0,\quad a=b=6,  \quad a'=b'=-2,\quad \chi(\EEE')=4.
\end{equation}
\end{lemma}

\begin{proof}
By \eqref{eq:c1c2} and \eqref{eq-dim-Fn:new-a} we have 
\begin{eqnarray*}
c'&=&c+eap-aq-bp-ep^2+2pq
\\
\nonumber
&\le& -\frac12ea(a+1)+ab+a+b+2- (g-1)+eap-aq-bp-ep^2+2pq
\\
\nonumber
&=& -\frac{1}{4} ea^2 -\frac{1}{2} ea -\frac{1}{4}ea'^2 + \frac{1}{2} ab +a +b + \frac{1}{2} a'b'-34.
\end{eqnarray*}
Recall that $b\le 6+3e$. Hence,
\begin{eqnarray}
\label{eq:new:cc}
c'&\le& -\frac{1}{4} ea^2 -\frac{1}{2} ea -\frac{1}{4}ea'^2 + \frac{1}{2} a(6+3e) +a +6+3e + \frac{1}{2} a'b'-34
\\
\nonumber
&=& -\frac{1}{4} ea^2 + (4 +e)a +3e -\frac{1}{4}ea'^2 + \frac{1}{2} a'b'-28.
\end{eqnarray}

Assume that $e=0$. Then $c'\le 4a + \frac{1}{2} a'b'-28\le -2$ because $a\le 6$. Moreover, 
the equality holds only if $a=6$, $a'=b'=-2$.
We obtain \eqref{new:lemma:c2:case1}.

Assume that $e=2$. Then 
\[
-2\le c'\le -\frac{1}{2} a^2 + 6a +6 -\frac{1}{2}a'^2 + \frac{1}{2} a'b'-28\le 24 -\frac{1}{2}a'^2 + \frac{1}{2} a'b'-28\le -\frac{5}2,
\]
a contradiction.

Finally, let $e=3$. Then 
\[
-2\le c' \le -\frac{3}{4} a^2 + 7a +9 -\frac{3}{4}a'^2 + \frac{1}{2} a'b'-28
\le - \frac{11}{4} -\frac{3}{4}a'^2 + \frac{1}{2} a'b'.
\]
Hence, $3a'^2 +3\le 2 a'b'$. This is impossible since $-2\le a',\, b'\le -1$.
\end{proof}

\begin{lemma}
$H^0(\EEE')\neq 0$.
\end{lemma}

\begin{proof}

By \eqref{eq-dim-Fn:new2} and Lemma~\ref{new:lemma:c2} we have
\[
\chi(\EEE')= \left(b'-\frac12ea'\right)(a'+1)+a'-c'+2>0.
\]
Assume that $H^0(\EEE')=0$. 
Then
$H^2(\EEE')\neq 0$. By Serre Duality
\[
0\neq H^2(\EEE')^{\vee}\simeq H^0(\EEE'^{\vee}\otimes \omega_Z)\simeq
H^0(\EEE'\otimes \det \EEE'^{\vee}\otimes \omega_Z).
\]
On the other hand, $(\det \EEE'^{\vee}\otimes
\omega_Z)^{\vee}=\OOO_Z((a'+2)\Sigma+(b'+(e+2))l)$, hence $H^0((\det
\EEE'^{\vee}\otimes \omega_Z)^{\vee})\neq 0$, a contradiction.
\end{proof}

\begin{proof}[Proof of Theorem~~\xref{th-main} under the assumptions~\ref{new-ass-0}--\ref{new-ass-3}]
Consider a nonzero section $s\in H^0(\EEE')$. If~$s$ does not
vanish anywhere, then $\EEE'$ is an extension of some line bundle
$\EEE_1$ by $\OOO$. But then $c_2(\EEE')=c_2(\EEE_1)=0$. This
contradicts Lemma~\ref{new:lemma:c2}. Then the zero locus of $s$ contains a curve $Y$,
because $\mathrm{c}_2(\EEE')<0$. Write $Y\sim
q_1 \Sigma+q_2 l$. Then by \eqref{eq:c1c2} restrictions $\EEE'$ to general curves
$l\in |l|$ and $\Lambda\in |\Sigma+el|$ have the form
\begin{eqnarray*}
\EEE'|_{l}&=&\OOO(q_1)\oplus\OOO(a'-q_1),
\\
\EEE'|_{\Lambda}&=&\OOO(q_2)\oplus\OOO(b'-q_2).
\end{eqnarray*}  
Taking Lemma~\ref{lemma-p1-p2-nachalo} into account we obtain\footnote{
In the journal version of this paper it was an error in the system of inequalities~\eqref{new-eq}. 
We would like to thank Haidong Liu and and Wenyou Li for pointing out this gap.} 
\begin{equation}
\label{new-eq}
\begin{array}{rcl}
2q_1-a'&\le &2,
\\
2q_2-b' &\le& 2+e.
\end{array}
\end{equation}
This immediately implies that
\[
e\ge 2,\qquad q_1=0, \qquad 1\le q_2\le 2. 
\]
In particular, $Y\sim q_2l$. 
Consider the vector bundle $\EEE'':=\EEE'\otimes \OOO_Z(-Y)$.
Then $H^0(Z,\EEE'')\neq 0$ and the zero locus of any section $s\in H^0(Z,\EEE'')$
does not contain curves. Therefore, 
$0\le \mathrm{c}_2(\EEE'')=\mathrm{c}_2(\EEE')-q_2a'$ and so 
\[
-3\ge \mathrm{c}_2(\EEE')\ge q_2a'
\]
by Lemma~\ref{new:lemma:c2}. We obtain $q_2=2$, $a'=-2$.
Now from \eqref{new-eq} we have  $e=3$ and $b'=-1$,
$c'\ge -4$.
In this case \eqref{eq:new:cc} has the form
\[ 
-4\le c'\le -\frac{3}{4} a^2 + 7a -21,
\]
which is impossible.

 Proposition
\xref{prop-surf-main} and Theorem~\xref{th-main} are proved.
\end{proof}

\end{document}